\documentclass[12pt]{amsart}
\usepackage{amsmath}
\usepackage{amsfonts}
\usepackage{amsthm}
\usepackage{amssymb}
\usepackage{amscd}
\usepackage[all]{xy}
\usepackage{enumerate}

\keywords{Waring problem, Variety of power sums, 
Fano threefold} 

\subjclass{Primary 14J45; Secondary 14N05, 14H42}

\theoremstyle{plain}
\newtheorem{thm}{Theorem}[subsection]

\newtheorem{prop}[thm]{Proposition}

\newtheorem{cor}[thm]{Corollary}
\newtheorem{lem}[thm]{Lemma}
\newtheorem{cla}[thm]{Claim}

\theoremstyle{definition}
\newtheorem{defn}[thm]{Definition}

\newtheorem{que}[thm]{Question}

\newtheorem{nota}[thm]{Notation}

\newtheorem*{ackn}{Acknowledgment}

\theoremstyle{remark}
\newtheorem*{rem}{Remark}
\newcommand{\sB}{\mathcal{B}}
\newcommand{\sC}{\mathcal{C}}
\newcommand{\sD}{\mathcal{D}}
\newcommand{\sE}{\mathcal{E}}
\newcommand{\sF}{\mathcal{F}}

\newcommand{\sH}{\mathcal{H}}
\newcommand{\sI}{\mathcal{I}}

\newcommand{\sL}{\mathcal{L}}
\newcommand{\sN}{\mathcal{N}}
\newcommand{\sM}{\mathcal{M}}
\newcommand{\sO}{\mathcal{O}}

\newcommand{\sQ}{\mathcal{Q}}

\newcommand{\sU}{\mathcal{U}}
\newcommand{\sV}{\mathcal{V}}

\newcommand{\sZ}{\mathcal{Z}}

\newcommand{\mC}{\mathbb{C}}
\newcommand{\mF}{\mathbb{F}}

\newcommand{\mN}{\mathbb{N}}
\newcommand{\mP}{\mathbb{P}}

\newcommand{\mZ}{\mathbb{Z}}

\newcommand{\Ima}{\mathrm{Im}\,}

\newcommand{\ap}{\mathrm{ap}}
\newcommand{\Bs}{\mathrm{Bs}\,}

\newcommand{\Aut}{\mathrm{Aut}\,}
\newcommand{\Hilb}{\mathrm{Hilb}}
\newcommand{\Hom}{\mathrm{Hom}}
\newcommand{\sHom}{\mathcal{H}{\!}om\,}

\newcommand{\Pic}{\mathrm{Pic}\,}

\newcommand{\Supp}{\mathrm{Supp}\,}

\newcommand{\VSP}{\mathrm{VSP}\,}

\numberwithin{equation}{section}

\title{On blow-ups of the quintic del Pezzo $3$-fold and \\
varieties of power sums of 
quartic hypersurfaces}

\author{Hiromichi Takagi}

\address{Graduate School of Mathematical Sciences \\
the University of Tokyo\\
Tokyo, 153-8914, Japan\\
\texttt{takagi@ms.u-tokyo.ac.jp}}

\author{Francesco Zucconi}

\address{D.I.M.I. \\
the University of Udine\\
Udine, 33100 Italy\\
\texttt{Francesco.Zucconi@dimi.uniud.it}}

\date{8.14, 2008}

\begin{document}

\maketitle

\markboth{Takagi and Zucconi}{Varieties of power sums}

\begin{abstract}
We construct new subvarieties in the varieties of power sums for
certain quartic hypersurfaces. This provides 
a generalization of Mukai's description of 
smooth prime Fano threefolds of genus twelve as
the varieties of power sums for plane quartics. 
In fact in \cite{TZ} we show that 
these quartics are exactly the Scorza quartics associated to general
pairs of trigonal curves and ineffective theta characteristics and
this enables us to prove there the main cojecture of \cite{DK}. 
\end{abstract}

\tableofcontents

%%%%%%%%%%%%%%%%%%%%%%%%%%%%%%%%%%%%%%%%%%%%%%%%%%%%%%%%%%%%%%%%%%%%%%%%%%%%%%%%%%%%%%%%%%%%%%%%%%%%%%%%%%%%%%%%%%%%%%%%%%%%%%%%%%%%%%%%%%%%%%%%%%%%%%%%%%%%%%%%
\section{Introduction}

\subsection{Varieties of power sums}
\label{subsection:intro1}~

The problem of representing a homogeneous form
as a sum of powers of linear forms 
has been studied since the last decades of the $19^{\rm{th}}$
century. This is called the {\em{Waring problem}} 
for a homogeneous form.
We are interested in the study of the global structure
of a suitable compactification of the variety 
parameterizing all such representations of
a homogeneous form.
A precise definition of the claimed compactification is the following:
\begin{defn}
Let $V$ be a $(v+1)$-dimensional vector space and 
let $F\in S^m \check{V}$ be a homogeneous forms of degree $m$ on $V$,
where $\check{V}$ is the dual vector space of $V$. Set
\[
\VSP(F,n)^o:=
\{([H_1],\dots, [H_n])\mid H_1^m+\cdots+H_n^m=F\}
\subset \Hilb^n (\mP_*\check{V}).
\]
The closed subset $\VSP(F,n):=\overline{\VSP(F,n)^o}$ is called  the {\em{varieties of power sums}} of $F$.
\end{defn}
Sometime $\mP_*\check{V}$ will be denoted by $\check{\mP}^v$. 

As far as we know,
the first global descriptions of positive dimensional
VSP's were given by Mukai.

\subsection{Mukai's result}
\label{subsection:Mukai}~

Let $A_{22}$ be a smooth prime Fano threefold of genus twelve,
namely, a smooth projective threefold
such that $-K_{A_{22}}$ is ample, the class of $-K_{A_{22}}$ generates $\Pic A_{22}$, 
and the genus $g(A_{22}):=\frac{(-K_{A_{22}})^3}{2}+1$ is equal to twelve.
The linear system $|-K_{A_{22}}|$ embeds $A_{22}$ into $\mP^{13}$.

Mukai discovered the following remarkable theorem (\cite{Mu2}, \cite{Mukai12}):
\begin{thm}\label{v22} Let $\{F_4=0\}\subset \mP^2$ 
    be a general plane quartic curve. Then 
    \begin{enumerate}[$(1)$]
\item
 ${\rm{VSP}}(F_4,6)\subset \Hilb^{6} \check{\mP}^{2}$
is a $A_{22}$; and conversely,
\item every general $A_{22}$ is of this form.
\end{enumerate}
\end{thm}
%Here the dual notation is used for later convenience.

His motivation to discover this result was
a characterization of a general $A_{22}$.
For this purpose, he noticed that
the Hilbert scheme of lines on a general $A_{22}\subset \mP^{13}$
is isomorphic to a smooth plane quartic curve $\sH_1\subset \mP^2$
(the notation $\mP^2$ will be compatible with $\check{\mP}^2$ in 
Theorem \ref{v22}).
He wanted to recover $A_{22}$ by $\sH_1$; for this, one more data was necessary.
In fact he proved that 
the correspondence on $\sH_1\times \sH_1$ defined 
by intersections of lines on $A_{22}$ gives
an ineffective theta characteristic $\theta$ on $\mathcal{H}_1$.
More precisely,
$\theta$ is constructed so that
the following two sets in $\sH_1\times \sH_1$ coincide:
\[
\{([l],[m])\mid l\cap m\not =\emptyset,l\not =m\}
=\{([l],[m])\mid h^0(\theta+[l]-[m])>0\}.
\]
Now a deep and beautiful result of G. Scorza asserts that,
associated to the pair $(\sH_1,\theta)$,
there exists another plane quartic curve $\{F_4=0\}$
in the same ambient plane as $\sH_1$.
(By saluting Scorza,
$\{F_4=0\}$ is called the {\em{Scorza quartic}}.)
Then, finally, Mukai proved that
$A_{22}$ is recovered as $\VSP(F_4, 6)$.
This is the result (2) of theorem \ref{v22}. We recall also that
since the number of the moduli of $A_{22}$ is equal to
$\dim \sM_4=6$, (1) follows from (2).

Moreover, Mukai observed that
conics on $A_{22}$ are parameterized by the plane $\sH_2$ and
$\sH_2$ is naturally considered as the plane $\check{\mP}^2$ 
dual to $\mP^2$ since, for a conic $q$ on $A_{22}$,
the lines intersecting $q$ form a hyperplane section of $\sH_1$.

Further, he showed that
the six points $[H_1],\dots, [H_6]$
such that $([H_1],\dots, [H_6])\in \VSP^o(F_4,6)$ 
correspond to six conics through one point of $A_{22}$. 

To sum up, even if it is not evident from the statement,
the content of Mukai's theorem is a new interpretation of 
the geometry of lines and conics on $A_{22}$.

\subsection{Generalization}
\label{subsection:main'}~

We study the relation between the concept of varieties of power sums 
and the geometry of lines and conics of other classes of $3$-folds. 

To do that, consider the smooth quintic del Pezzo threefold $B$
namely, a smooth projective threefold
such that $-K_{B}=2H$, where $H$ is the ample generator
of $\Pic B$ and $H^3=5$. It is well known that the linear system 
$|H|$ embeds $B$ into $\mP^6$.

Now, following Iskovskih we doubly project 
$A_{22}$ from a general line, that is we consider the following diagram:
\begin{equation*}
\xymatrix{
& A' \ar[dl]_{f'}  &\dashrightarrow & 
 A \ar[dr]^{f} & \\
 A_{22}  &  &  & & B, }
\end{equation*}
where
\begin{itemize}
\item
$f'$ is the blow-up along a general line $l$,
\item $A'\dashrightarrow A$ is a flop,
\item
$f$ is the blow-up along a smooth rational curve of degree five,
where the degree is measured by $H$. We consider $B\subset \mP^6$ by $\Phi_{|H|}$.
\end{itemize}
(See also the section \ref{section:v22} for more information).

It is known that a general line on $A_{22}$ is mapped to a general line on $B$
intersecting $C$, and
a general conic on $A_{22}$ is mapped to a general conic on $B$
intersecting $C$ twice.
These facts are easy to see since the exceptional divisor of $f$ is
the strict transform of the unique hyperplane section 
vanishing along $l$ with multiplicity $3$.
  
This situation is generalizable by considering
a general smooth rational curve $C$ of degree $d$ on $B$, 
where $d$ is an arbitrary integer greater than or equal to $5$
(mainly $d\geq 6$) and the sets of the
secant lines of $C$ and of the multi-secant conics of $C$ respectively.
This led to the following definition:

\begin{defn}
    \begin{enumerate}[(1)]
\item
A pair $(l,t)$ of a line $l$ on $B$
and a point $t\in C\cap l$ is called a {\em{marked line}}.
\item
A pair of a conic $q$ on $B$ 
and a zero-dimensional subscheme $\eta \subset C$
of length two contained in $q_{|C}$ 
is called a {\it{marked conic}}.
\end{enumerate}
\end{defn}

We can prove:

\begin{prop}
\label{prop:intro1}
Marked lines are parameterized by
a smooth trigonal canonical curve
$\mathcal{H}_1$ of genus $d-2$.
\end{prop}

See the subsection \ref{subsection:H1} for the proof.
Here is a sketch of the proof.
It is known that there are three lines (counted with multiplicities)
through a point of $B$ (see the subsection \ref{subsection:FN}).
This gives the triple cover $\sH_1\to C$
such that $(l,t)\mapsto t$.
Moreover, points where `special lines' pass through form a divisor $\in |2H|$
and the intersection of this divisor and $C$ is nothing but 
the branch locus of this triple cover. 
We can show that all ramifications are simple.  
Thus it holds
\[
2g(\sH_1)-2=3(-2)+2d,\ \text{namely,} \ g(\sH_1)=d-2.
\]

As Mukai did, we can
define an ineffective theta characteristic $\theta$ on $\sH_1$
and construct the Scorza quartic hypersurface $\{F_4=0\}$
associated to this in the sense of
\cite[\S 9]{DK}.
This quartic hypersurface lives 
in the projective space $\mP^{d-3}\supset \sH_1$.
This construction, however, is rather indirect, hence
we give a more direct construction of $F_4$ in this paper. 
We will show the quartic constructed in this paper is actually Scorza
in the forthcoming paper \cite{TZ}. 

For the construction of the quartic $\{F_4=0\}$,
we make use of marked conics, which we study in the subsection 
\ref{subsection:H2} in detail. Among other things, we prove the following:

\begin{prop}
If $d\geq 6$, then marked conics are parameterized by
a so-called {\em{White surface}} $\sH_2$ 
obtained by blowing up $S^2 C\simeq \mP^2$ at $\binom{d-2}{2}$ points.
$\sH_2$ is embedded by $|(d-3)h-\sum_{i=1}^{s} e_i|$ into $\check{\mP}^{d-3}$,
where $h$ is the pull-back of a line, $e_i$ 
are the exceptional curves of $\sH_2\to \mP^2$
and $s:=\binom{d-2}{2}$. 
\end{prop}

Here we use the notation $\check{\mP}^{d-3}$
since the ambient projective spaces of $\sH_1$ and $\sH_2$
are reciprocally dual as in the Mukai's case.
If $d=6$, then $\sH_2$ is a cubic surface.
In general,
Gimigliano \cite{Gimi} shows that $\sH_2$ is the intersection of cubics.

The proof of this proposition
is more involved than that of Proposition \ref{prop:intro1}.
See Corollary \ref{cor:H2} and Theorem \ref{thm:H_2} for the proof.
Here is a sketch of the proof.
The morphism $\sH_2\to \mP^2$ is just a 
natural one $\sH_2\to S^2 C\simeq \mP^2$
mapping $(q,\eta)\mapsto \eta$.
Let $\beta_i$ be a bi-secant line of $C$.
It is shown that there exist
$s:=\binom{d-2}{2}$ bi-secant lines of $C$
(see Corollary \ref{orabisecanti}).
Then for the length two subscheme $\beta_{i|C}$,
there exist infinitely many marked conics
$(\beta_i\cup \alpha, \beta_{i|C})$,
where $\alpha$ are lines intersecting $\beta_i$, and
it is known that such $\alpha$'s form one-dimensional family
(see Proposition \ref{prop:FN} (5)).
This indicates why $\sH_2\to S^2 C$ is the blow-up at 
$s$ points, which are $[\beta_{i|C}]\in S^2 C$. 
Moreover, birationality of $\sH_2\to \mP^2$ follows from the fact
that there exists a unique conic on $B$ through two points $t_1$ and $t_2$
if there is no line on $B$ through $t_1$ and $t_2$.
This can be seen by the double projection from $t_1$ 
(see Corollary \ref{cor:uniqueconic}).

Actually we consider the curves on $A$
called lines and conics on $A$ corresponding one to one
to marked lines and conics respectively.

In \cite[\S 9]{DK},
the quartic $F_4$ is constructed for 
$(\sH_1,\theta)$, which is a data of intersections of marked lines.
Here to construct $F_4$ we need data of intersections of marked conics. 

In fact assume that $d\geq 6$.
Consider the locus
$D_l\subset \sH_2$ parameterizing marked conics 
which intersect a fixed marked line $l$.
The locus $D_l$ turns out to be a divisor
linearly equivalent to
$(d-3)h-\sum_{i=1}^s e_i$ on $\sH_2$.
Moreover,
$|D_l|$ is very ample and embeds $\sH_2$ in $\check{\mP}^{d-3}$
(see Theorem \ref{thm:H_2} (1)).
Set
$\sD_2:=\{([q_1],[q_2])\in\sH_{2}\times\sH_{2}\mid q_1\cap q_2\neq \emptyset\}$
and denote by $D_q$ the fiber of $\sD_2\to \sH_2$ over a point $[q]$.
It is easy to verify $D_q\sim 2D_l=\sO_{\sH_2}(2)$.
By the seesaw theorem,
it holds that $\sD_2\sim p_1^*D_q+p_2^*D_q$.
Since $\sH_2$ is projectively Cohen-Macaulay and 
is not contained in a quadric (Theorem \ref{thm:H_2} (4)),
it holds
$H^0(\sH_2\times \sH_2,\sD_2)\simeq 
H^0(\check{\mP}^{d-3}\times \check{\mP}^{d-3},
\sO(2,2))$.
Thus
$\sD_2$ is the restriction of a unique $(2,2)$-divisor $\sD'_2$ on 
$\check{\mP}^{d-3}\times \check{\mP}^{d-3}$.
Since $\sD'_2$ is symmetric, we may assume 
its equation $\widetilde{\sD}_2$ is also symmetric.
By restricting $\widetilde{\sD}_2$ to the diagonal, we obtain 
a quartic hypersurface $\{\check{F}_4=0\}$ in $\check{\mP}^{d-3}$.
We can show that $\check{F}_4$ is non-degenerate 
in the sense of \cite{doldual} (see the appendix).
Then there exists a unique quartic hypersurface $\{F_4=0\}$ 
in $\mP^{d-3}$ called the quartic form dual to $\check{F}_4$.

Now we can state our main result, which generalizes (2)
of Theorem \ref{v22}:
\begin{thm}
\label{thm:main2}
Let $f\colon A\to B$ be the blow-up along $C$, and
let $\rho\colon \widetilde{A}\to A$ be the blow-up
of $A$ along the strict transforms $\beta'_i$
of $\binom{d-2}{2}$ bi-secant lines $\beta_i$ of $C$ on $B$.
Then there is an injection from $\widetilde{A}$ to 
$\VSP(F_4,n)$,
where $n:=\binom{d-1}{2}$.
Moreover the image is uniquely determined by the incident variety $\sD_2$
and is an irreducible component of
\[
\VSP(F_4, n;\sH_2):=
\overline{
\{([H_1],\dots, [H_n])\mid [H_i] \in \sH_2\}}
\subset \VSP(F_4,n).
\]
\end{thm}
See Theorem \ref{diretto} and Proposition \ref{lem:main}.

Actually, the number $n$ is equal to the number of multi-secant conics of 
$C$ through a general point of $B$ (see Corollary \ref{cor:n}).
Moreover, rather importantly,
\begin{equation}
\label{eq:intro0}
\text{
$n$ is equal to the dimension of
quadric forms on $\check{\mP}^{d-3}$.}
\end{equation}

We give an outline of the proof of the main result.
Let $\sU_2\to \sH_2$ be the universal family of conics on $A$,
and consider the natural projection 
$\psi\colon \sU_2\subset A\times \sH_2 \to A$.
The morphism $\psi$ is not finite 
(see Proposition \ref{fuoripiatto}). Nevertheless
the blow-up $\widetilde{\sU}_2\to 
\sU_2$ along $(\cup \beta'_i\times \sH_2)\cap \sU_2$ is Cohen-Macaulay and
the natural projection 
$\widetilde{\psi}\colon \widetilde{\sU}_2\to \widetilde{A}$ is finite
of degree $n$ (Proposition \ref{finitezzafinale}).
Therefore, since $\widetilde{\sU}_2\subset \widetilde{A}\times \sH_2$,
$\widetilde{\psi}$ is a flat family of 
$0$-dimensional subschemes $\subset \sH_2$ of length $n$
parameterized by $\widetilde{A}$.
Geometrically, the fiber over a general point $\widetilde{a}\in \widetilde{A}$
corresponds to $n$ conics through the image of $\widetilde{a}$ on $A$.
The morphism $\widetilde{\psi}$ 
defines $\widetilde{A}\to \Hilb^n \check{\mP}^{d-3}$ which is the one 
claimed in the main theorem.
To understand its image, we need to understand the double polars of the 
special quartic $F_4$.

By the construction of $\check{F}_4$ and the theory of polarity
(see the appendix), it holds that,
for a conic $q$ on $A$ and the hyperplane section 
$\{H_q=0\}\subset \mP^{d-3}$ corresponding to the point 
$[q]\in \check{\mP}^{d-3}$, 
the locus $D_q$ is equal to
$\{\widetilde{D}_q:=P_{H_q^2} (\check{F}_4)=0\}\cap \sH_2$. 
By definition of the dual quartic form $F_4$,
it holds 
\begin{equation}
\label{eq:intro1}
P_{\widetilde{D}_q}(F_4)=H_q^2.
\end{equation}
Moreover,
by definition of $D_q$, it holds
that, for $n$ conics $q_1,\dots,q_n$ on $A$ corresponding
to a general fiber of $\widetilde{\psi}$,
\begin{equation}
\label{eq:intro2}
\text{$\widetilde{D}_{q_i} ([q_i])\not =0$ and
$\widetilde{D}_{q_i} ([q_j])=0$ $(i\not =j)$}.
\end{equation}
Now the main theorem follows from a more or less formal
argument of the theory of polarity from 
(\ref{eq:intro0}), (\ref{eq:intro1}),  and (\ref{eq:intro2}).

We believe that, even by reading the proof
of Theorem \ref{diretto} after 
reading only this introduction and possibly
the appendix,
the readers can understand 
at least the reason why the variety of power sums appears.

\subsection{Structure of the paper}~

We add some explanations about the structure of the paper.

In the section 2,
we construct smooth rational curves $C_d$ of degree $d$ on $B$
and study in detail the relation of general $C_d$ with lines and conics on $B$.

In the section 3, we describe the projection of $B$
from a line or a conic, and the double projection of $B$ from a point.
These operations are useful for counting the number of multi-secant conics of $C$
satisfying various pre-specified geometric conditions. 
For example, using double projection from a general point of $B$,
we can show that the number of multi-secant conics of $C$ through a 
general point of $B$ is
equal to $n$ (see Corllary \ref{cor:n}).

Sections 2 and 3 are rather technical as far as the proofs it concerns
but the results are really easy to be understood by a general reader
and at least one of them, we mean Proposition \ref{prop:Cd2}, is
of unexpected geometrical content; 
Proposition \ref{prop:Cd2} or its restatement 
Corollary \ref{cor:finite} shows that the number of multi-secant conics of $C$ through
{\em{any}} point of $B$ outside $C$ is finite.
This will be refined to finiteness results
contained into Propositions \ref{fuoripiatto} and \ref{finitezzafinale}.

In the section 4,
we mostly study marked lines and conics,
and lines and conics on 
the blow-up $A$ of $B$ along a smooth rational curve $C$ of degree $d$
as we mentioned in the subsection \ref{subsection:main'}. 

In the section 5, we show the main theorem.

In the section \ref{section:v22},
we explain Mukai's result from our view point.

Finally we add an appendix which forms the section 7,
where we explain some very basic facts on the theory of polarity 
for the readers' convenience.

\subsection{Forthcoming paper}~

This work laids the foundations for the results of
\cite{TZ}.

As we mentioned in the abstract,
there we show that the quartic $\{F_4=0\}$ 
coincides with the Scorza quartic
associated to $(\sH_1,\theta)$ and the theta charasterisctic
$\theta$ is constructed explicitely.

Following \cite{DK},
we also study other geometric objects associated to 
$(\sH_1,\theta)$.
As an amazing application,
we  show the existence of the Scorza quartics for any general pairs of
curves and ineffective theta characteristics. 
This is an affirmative answer to the conjecture stated by Dolgachev and Kanev in \cite[\S 9]{DK}.

Moreover, we can study the moduli spaces of spin curves, especially of
trigonal spin curves
relating this with the Hilbert schemes of smooth rational curves on $B$.
In fact we prove that $\sH_1$ is a general trigonal curve if $C$ is general.
\begin{ackn}
We are thankful to Professor S. Mukai for 
valuable discussions and constant interest on this paper.
We received various useful comments from 
K. Takeuchi, A. Ohbuchi, S. Kondo,
to whom we are grateful.
The first author worked on this paper
partially when he was staying at the Johns Hopkins University under the
program of Japan-U.S. Mathematics Institute (JAMI)
in November 2005 and at the Max-Planck-Institut f\"ur Mathematik
from April, 2007 until March, 2008. 
The authors 
worked jointly during the first author's stay at the Universit\`a di Udine
on August 2005, and the Levico Terme conference on Algebraic Geometry in
Higher dimensions on June 2007.
The authors are thankful to all the above institutes
for the warm hospitality they received. 
\end{ackn}
%
%%%%%%%%%%%%%%%%%%%%%%%%%%%%%%%%%%%%%%%%%%%%%%%%%%%%%%%%%%%%%%%%%%%%%%%
%%%%%%%%%%%%%%%%%%%%%%%%%%%%%%%%%%%%%%%%%%%%%%%%%%%%%%%%%%%%%%%%%%%%%%%%%%%%%%%%%%%%%%%%%%%%%%%%%%%%%%%%%%%%%%%%%%%%%%%%%%%%%%%%%%%%%%%%%%%%%%%%%%%%%%%%%%%%%%%%%%%%%%%%%%%%%%%%%%%%%%%%%%%%%%%%%%%%%%%%%%%%%%%%%%%%%%%%%%%%%%%%%%%%%%%%%%%%%%%%%%%%%%%%%%%%%%%%%%%%%%%%%%%%%%%%%%%%%%%%%%%%%%%%%%%%%%%%%%%%%%%%%%%%%%%%%%%%%%%%
\section{Rational curves on the quintic del Pezzo threefold $B$}

Let $V$ be a vector space with ${\rm{dim}}_{\mC}V=5$. The
Grassmannian $G(2,V)$ embeds into $\mP^{9}$ and we denote the image 
by $G\subset\mP^{9}$. It is well-known that the quintic del Pezzo $3$-fold,
i.e.,
the Fano $3$-fold $B$ of index $2$ and of degree $5$ 
can be realized as $B=G\cap \mP^{6}$, where $\mP^{6}\subset \mP^9$ 
is transversal to $G$ (see \cite{Fu2}, 
\cite[Thm 4.2 (iii), the proof p.511-p.514]{I1}).

First we collect basic known facts on lines and conics on $B$ 
almost without proof.
Let $\sH_{{1}}^{B}$ and $\sH_{{2}}^{B}$ be the Hilbert scheme, 
respectively, of lines and of conics on $B$.

\subsection{Lines on $B$}
\label{subsection:FN}~

Let $\pi\colon \mP\to \sH_{1}^{B}$ be the universal
family of lines on $B$ and 
$\varphi\colon \mP\to B$ the natural projection. 
By \cite[Lemma 2.3 and Theorem I]{FuNa}, $\mathcal{H}_{1}^{B}$ is isomorphic to $\mP^2$
and $\varphi$ is a finite morphism of degree three. 
In particular the number of lines passing through a point 
is three counted with multiplicities.  
We recall some basic facts about $\pi$ and $\varphi$ which we use in
the sequel.

Before that, we fix some notation.
\begin{nota}
For an irreducible curve $C$ on $B$, 
denote by $M(C)$ the locus $\subset \mP^2$ of lines intersecting $C$,
namely, $M(C):=\pi(\varphi^{-1}(C))$ with reduced structure. 
Since $\varphi$ is flat, $\varphi^{-1}(C)$ is purely one-dimensional.
If $\deg C\geq 2$, 
then $\varphi^{-1}(C)$ does not contain a fiber of 
$\pi$, thus $M(C)$ is a curve.
See Proposition \ref{prop:FN} for the description of $M(C)$
in case $C$ is a line.  
\end{nota}

\begin{defn}
A line $l$ on $B$ is called a {\em{special line}} if
$\sN_{l/B}\simeq 
\sO_{\mP^1}(-1)\oplus \sO_{\mP^1}(1).$
\end{defn}

\begin{rem}
If $l$ is not a special line on $B$,
then 
$\sN_{l/B}=\sO_{l}\oplus\sO_{l}$.
\end{rem}

\begin{prop} 
\label{prop:FN}
It holds$:$ 
\begin{enumerate}[$(1)$]
 \item
 for the branched locus $B_{\varphi}$ of 
 $\varphi\colon \mP\to B$ we have:
\begin{enumerate}[$({1}\text{-}1)$]
\item $B_{\varphi}\in |-K_{B}|$, and
\item $\varphi^*B_{\varphi}=R_1+2R_2$,
where $R_1\simeq R_2\simeq \mP^1\times \mP^1$, and
$\varphi\colon R_1\to
B_{\varphi}$ and 
$\varphi\colon R_2\to
B_{\varphi}$ are injective,
\end{enumerate}
\item $R_{2}$ is contracted to a conic $Q_2$ by 
$\pi\colon \mP\to \sH_{1}^{B}$. Moreover $Q_{2}$ is 
the branched locus of the finite double cover 
$\pi_{|R_1}\colon R_{1}\rightarrow
\mathcal{H}^{B}_{1}$,
\item $Q_2$ parameterizes special lines,
\item
if $l$ is a special line,
then $M(l)$ is the tangent line to $Q_2$ at $[l]$.
If $l$ is not a special line,
then $\varphi^{-1}(l)$ is the disjoint union of 
the fiber of $\pi$ corresponding to $l$,
and the smooth rational curve dominating a line on $\mP^2$.
In particular,
$M(l)$ is the disjoint union of a line and the point $[l]$.

{\em{By abuse of notation}}, we denote by $M(l)$ the one-dimensional part of
$M(l)$ for any line $l$. 
Vice-versa,
any line in $\mathcal{H}_{1}^{B}$ is of the form $M(l)$ for some line $l$,
and
\item 
the locus swept by lines intersecting $l$ is a hyperplane section $T_{l}$
  of $B$ whose singular locus is $l$. For every point $b$ of $T_l\setminus l$,
  there exists exactly one line which belongs to $M(l)$ 
  and passes through $b$.
Moreover, if $l$ is not special, then the normalization of $T_l$ is $\mF_1$
and the inverse image of the singular locus is the negative section of $\mF_1$,
or, if $l$ is special, then 
the normalization of $T_l$ is $\mF_3$
and the inverse image of the singular locus is the union of 
the negative section and a fiber.
\end{enumerate}
\end{prop}
\begin{proof} See \cite[\S 2]{FuNa} and \cite[\S 1]{IlievB5}. 
\end{proof}

By the proof of \cite{FuNa} we see that $B$ is stratified
according to the ramification of $\varphi\colon \mP\to B$ as follows:
\[
B=(B\setminus B_{\varphi})\cup 
(B_{\varphi}\setminus C_{\varphi})\cup C_{\varphi},
\] where 
$C_{\varphi}$ is a smooth rational normal sextic and 
if $b\in B\setminus B_{\varphi}$ exactly three
distinct lines pass through it, if $b\in
(B_{\varphi}\setminus 
C_{\varphi})$ exactly two distinct lines pass through it,
one of them is special,
and finally $C_{\varphi}$ is the loci of $b\in B$ 
through which it passes only one line, which is special.

\subsection{Conics on $B$}

\begin{prop}
\label{prop:conics}
The Hilbert scheme $\sH^B_2$ of conics on $B$ is isomorphic to 
$\mP^4=\mP_*\check{V}$. The support of a double line is a special line and
the double lines are parameterized by a
rational normal quartic curve $\Gamma\subset \mP_*\check{V}$ and the
secant variety of $\Gamma$ is a singular cubic hypersurface which
is the closure of the loci parameterizing reducible conics.
\end{prop}
\begin{proof} See \cite[Proposition 1.2.2]{IlievB5}.

The identification in the first statement is given by the map 
$sp\colon \sH_{2}^{B}\rightarrow \mP_*\check{V}$ with $[c]
\mapsto\langle Gr(c)\rangle=\mP^{3}_{c}\subset\mP_*V$, where for a
general conic $c\subset B$ we set
\[
Gr(c):=\cup\{r\in \mP_*V\mid [r]\in c\}\simeq \mP^{1}\times\mP^{1}.
\]  

\end{proof}

\subsection{Construction of rational curves $C_d$ of 
degree $d$ on $B$}
\label{subsection:CdB}
~

We construct smooth rational curves of degree $d$ on $B$ 
by smoothing the union of a smooth rational curve of degree $d-1$
and one of its uni-secant lines.

\begin{defn}
Let $C$ and $\gamma$ be smooth curves on $B$.
We say that
$\gamma$ is a secant curve of $C$ if $C\cap \gamma\not =\emptyset$.
Moreover, we say that
$\gamma$ is a $k$-secant curve (resp. a multi-secant curve)
if $\gamma_{|C}$ is a $0$-dimensional subscheme of length $k$
(resp. of length greater than or equal to $2$). 
For $k=1,2,\dots$, we say uni-secant, bi-secant, \dots, instead.
\end{defn} 

\begin{prop}
\label{prop:Cd0}
There exists a smooth rational curve $C_d$ of degree $d$ on $B$ such that
\renewcommand{\labelenumi}{\textup{(\alph{enumi})}}
\begin{enumerate}
\item
a general line on $B$ intersecting $C_d$ is uni-secant,
\item 
$C_d$ is obtained as a smoothing of 
the union of a smooth rational curve $C_{d-1}$ of degree $d-1$
on $B$ and a general uni-secant line of it on $B$, and
\item
$\sN_{C_{d}/B}\simeq\sO_{\mP^{1}}(d-1)\oplus\sO_{\mP^{1}}(d-1)$. 
In particular 
$h^1(\sN_{C_{d}/B})=0$ and $h^0(\sN_{C_{d}/B})=2d$.
\end{enumerate}
\end{prop}
\begin{proof}
    We argue by induction on $d$.

If $d=1$, we have the assertion since 
$\sN_{C_1/B}\simeq\sO_{\mP^{1}}\oplus\sO_{\mP^{1}}$
for a general line $C_1$.

Now assume that $C_{d-1}$ is a smooth rational curve of degree $d-1$
on $B$ constructed inductively.
By induction, a general secant line $l$ of $C_{d-1}$ on $B$ is uni-secant.
Set $Z:=C_{d-1}\cup l$ and $\sN_{Z/B}: =\sHom_{\sO_{B}}(\sI_{Z},\sO_{B})$.
By induction, the normal bundle of $C_{d-1}$ satisfies (c). 
Thus, by $\sN_{l/B}\simeq\sO_{\mP^{1}}\oplus\sO_{\mP^{1}}$
and \cite[Theorem 4.1 and its proof]{HH},
it holds $h^1(\sN_{Z/B})=0$, and moreover
$Z:=C_{d-1}\cup l$ is strongly smoothable,
namely, we can find a smoothing $C_d$ 
of $Z$ with the smooth total space.
By the upper semi-continuity theorem,
$h^1(\sN_{C_d/B})=0$ and, by the Riemann-Roch
theorem, $h^0(\sN_{C_d/B})=2d$.

We check the form of the normal bundle of $C_d$. 
Set 
$\sN_{C_{d}/B}:=\sO_{\mP^{1}}(a_{d})\oplus\sO_{\mP^{1}}(b_{d})$ ($a_d\geq b_d$)
for the smoothing $C_d$ of $Z$. 
We show that $a_d=b_d=d-1$.
It suffices to prove $h^0(\sN_{Z/B}(-d))=0$.
In fact, then, by the upper semi-continuous theorem,
we have $h^0(\sN_{C_d/B}(-d))=0$ and $a_d, b_d\leq d-1$.
Thus, 
by $a_d+b_d=2d-2$, it holds $a_d=b_d=d-1$.
By noting
$\sN_{C_{d-1}/B}=\sO_{\mP^{1}}(d-2)\oplus\sO_{\mP^{1}}(d-2)$,
the equality $h^0(\sN_{Z/B}(-d))=0$ easily 
follows from the following three exact sequences, 
where $t:=C_{d-1}\cap l$:

\begin{equation*}
    0\rightarrow \sN_{Z/B}\rightarrow
    \sN_{Z/B\vert C_{d-1} }\oplus\sN_{Z/B\vert l}\rightarrow
    \sN_{Z/B}\otimes_{\sO_{B}}\sO_t\rightarrow 0.
 \end{equation*}
 
 \begin{equation*}
    0\rightarrow \sN_{C_{d-1}/B}\rightarrow
    \sN_{Z/B\vert C_{d-1}}\rightarrow T^{1}_t\rightarrow 0.
 \end{equation*}

 \begin{equation*}
    0\rightarrow \sN_{l/B}\rightarrow
    \sN_{Z/B\vert l}\rightarrow T^{1}_t\rightarrow 0.
 \end{equation*}
 
 We can inductively show that a general line $m$ intersecting $C_{d-1}$
does not intersect $l$, thus $m$ is a uni-secant line of $C_{d-1}\cup l$.
This implies (a) for $C_d$ by a deformation theoretic argument.

\end{proof}

\begin{cor}
\label{cor:smooth}
Let $C_d$ be a smooth rational curve of degree $d$
constructed as in Proposition $\ref{prop:Cd0}$.
The Hilbert scheme of smooth rational curves on $B$ of degree $d$
is smooth at $[C_d]$ and is of dimension $2d$. 
\end{cor}

\begin{proof}
The assertion follows from Proposition \ref{prop:Cd0} (c).
\end{proof}

\subsection{Relations of $C_d$ with lines and conics}

We study multi-secant lines and conics of $C_d$.

\begin{prop}
\label{prop:Cd1}
A general $C_d$ as in Proposition $\ref{prop:Cd0}$
satisfies the following conditions$:$
\begin{enumerate}[$(1)$]
\item
there exist no $k$-secant lines of $C_d$ on $B$ with $k\geq 3$,
\item there exist at most finitely many bi-secant lines of $C_{d}$ on
$B$, and any of them intersects $C_{d}$ simply,
\item bi-secant lines of $C_{d}$ on
$B$ are mutually disjoint,
\item neither a bi-secant line nor
a line through the intersection point between a bi-secant line and $C_d$
is a special line, and
\item
there exist at most finitely many points $b$ outside $C_d$
such that all the lines through $b$ intersect $C_d$,
and such points exist outside bi-secant lines of $C_d$.
\end{enumerate}
\end{prop}

\begin{proof}
We can prove the assertions 
by simple dimension counts based upon Proposition \ref{prop:Cd0}.
We assume that $d\geq 4$ since otherwise
we can verify the assertion easily.

$(1)$. Let $\sD$ be the closure of the set
\begin{align*}
\{([C_d], [l])\mid  C_d\cap l\ \text{consists of
$3$ points}
\}  \subset \sH_d^{B}\times \sH_1^{B}.
\end{align*}
Let $\pi_{d}\colon \sD\rightarrow\sH_d^{B}$ and  
$\pi_{1}\colon \sD\rightarrow\sH_1^{B}$ be the natural
morphisms induced by
the projections. The claim follows if we show that
${\rm{dim}}_{\mC}\sD\leq 2d-1$
since $\dim \sH_d^{B}=2d$.

Thus we estimate $\dim_{\mC} \Hom^{2d} (\mP^1, B;p_{i}\mapsto s_{i}, i=1,2,3)$ 
at $[\pi]$, 
where $p_{i}$, $i=1,2,3$ are fixed  points of $\mP^1$, $[\pi]$ 
is a general point and the degree is measured by $-K_{B}$. 
By $d\geq 4$ and Proposition \ref{prop:Cd0} (c), it holds that 
$h^{0}(\mP^{1},\pi^*T_{B}(-p_{1}-p_{2}-p_{3}))=2d-6$ and
$h^{1}(\mP^{1},\pi^*T_{B}(-p_{1}-p_{2}-p_{3}))=0$. 
Then 
\[
\dim_{\mC} \Hom^{2d} (\mP^1, B, p_{i}\mapsto s_{i}, i=1,2,3)_{[\pi]}
=h^0(\pi^*T_{B} (-p_{1}-p_{2}-p_{3}))
=2d-6.
\]
This implies that $\dim_{\mC}\pi_{1}^{-1}([l])\leq 2d-6 +3=2d-3$ since
the three points can be chosen arbitrarily. Then ${\rm{dim}}_{\mC}\sD\leq 2d-1$
since ${\rm{dim}}_{\mC}\sH_1^{B}=2$.
    
$(2)$. Now let 
$\sD$ be the closure of the set
\begin{align*}
\{([C_d], [l])\mid C_d\cap l\ \text{consists of
$2$ points} \}  
\subset \sH_d^{B}\times \sH_1^{B}.
\end{align*}
As before, let $\pi_{d}\colon \sD\rightarrow\sH_d^{B}$ and  
$\pi_{1}\colon \sD\rightarrow\sH_1^{B}$ be the natural
morphisms induced by
the projections.
By $d\geq 4$ and Proposition \ref{prop:Cd0} (c), it holds that
$h^{0}(\mP^{1},\pi^*T_{B}(-p_{1}-p_{2}))=2d-3$ and 
$h^{1}(\mP^{1},\pi^*T_{B}(-p_{1}-p_{2}))=0$. Then
\[
\dim_{\mC} \Hom^{2d} (\mP^1, B, p_{i}\mapsto s_{i}, i=1,2)_{[\pi]}=
h^0(\pi^*T_{B} (-p_{1}-p_{2}))
=2d-3.
\]
Since 
$\dim_{\mC}\Aut (\mP^1, p_{1},p_{2})=1$ it holds that
$\dim_{\mC}\pi_{1}^{-1}([l])\leq 2d-3+2-1=2d-2$. Hence 
${\rm{dim}}_{\mC}\sD=2d$. 
Hence $C_{d}$ has only a
finite number of bi-secant lines. 

We now show that the loci where $C_{d}$
has a tangent bi-secant is a codimension one loci inside $\sH_d^{B}$.
Let $B_{t}$ be the blow-up of $B$ in a point $t\in C_{d}$ and let $l$ 
be a bi-secant which is tangent to $C_{d}$ at $t$ (if it exists). Let $E$ be 
the exceptional divisor, and $C'$ and $l'$ the strict transforms of $C$ and 
$l$ respectively. By hypothesis there exists a unique point $s\in
E\cap C'\cap l'$. We estimate $\dim_{\mC} \Hom^{d-2} (\mP^1, B_{t}, 
p\mapsto s)_{[\pi]}$, where
$p$ is a fixed  point of $\mP^1$, $[\pi]$ is a general point, and 
the degree is measured by $-K_{B_t}$. 
In this case $h^0(\pi^*T_{B_{t}} (-p))= 2d-2$
hence $\dim_{\mC}\pi_{1}^{-1}([l])\leq 2d-2+ 1 -2=2d-3$. This implies 
the claim.

The cases $(3)$, $(4)$ and $(5)$  are similar. 
Thus we only give few comments for (5).
Set $\sD$ be the closure of the set
\begin{align*}
\{([C_d], [l_1],[l_2],[l_3])\mid 
C_d\cap l_i\not =\emptyset\, (i=1,2,3),\\l_1\cap l_2\cap l_3\not =\emptyset,
l_1\cap l_2\cap l_3\not \in C_d,
l_i \ \text{are distinct}\} \\
\subset \sH_d^{B}\times \sH_1^{B}\times \sH_1^{B}\times \sH_1^{B}.
\end{align*}
For the former half of (5), we have only to prove that $\dim \sD\leq 2d$.
This can be carried out by a similar dimension count as above.
For the latter half of (5), we use the inductive construction 
of $C_d$ besides dimension count. 
\end{proof}

We can prove the following by a similar method hence
we omit the proof.

\begin{prop}
\label{prop:Cd1'}
A general $C_d$ as in Proposition $\ref{prop:Cd0}$
satisfies the following conditions$:$
\begin{enumerate}[$(1)$]
\item 
there exist no $k$-secant conics of $C_{d}$ with $k\geq 5$,
\item
there exist at most finitely many quadri-secant conics of $C_d$ on
$B$, and no quadri-secant conic is tangent to $C_d$, and
\item $q_{|C_d}$ has no point of multiplicity greater than two
for any multi-secant conic $q$.
\end{enumerate}
\end{prop}

\begin{nota}
The bisecant lines of $C_{d}$ are denoted by 
$\beta_i$ where $i=1,\cdots, s$.
\end{nota}

In the following proposition,
we describe some more relations of $C_d$ with lines on $B$
which can be translated into the geometry of $\sH^B_1$.
More explicitly, we prove that $M(C_d)$ is suffciently general
if $C_d$ is general
(recall the notation of the subsection \ref{subsection:FN}).
\begin{prop}
\label{prop:Cd}
A general $C_d$ as in Proposition $\ref{prop:Cd0}$
satisfies the following conditions$:$

\begin{enumerate}[$(1)$]
\item
$C_d$ intersects $B_{\varphi}$ simply,

\item 
$M_d:=M(C_d)$ intersects $Q_2$ simply,
\item
$M_d$ is an irreducible curve of degree $d$
with only simple nodes
$($recall that in Proposition $\ref{prop:FN}$ $(4)$,
we abuse the notation by denoting the one-dimensional part of
$\pi(\varphi^{-1}(C_1))$ by
$M(C_1)$$)$,
\item
for a general line $l$ intersecting $C_d$,
$M_d\cup M(l)$ has only simple nodes as its singularities, and
\item
$M_d\cup M(\beta_i)$ has only simple nodes as its singularities.

\end{enumerate}
\end{prop}

\begin{proof}
We show the assertion inductively using
the smoothing construction of $C_d$ from
the union of $C_{d-1}$ and a general uni-secant line $l$
of $C_{d-1}$.

In case of $d=1$, by letting $C_1$ be a general line,
the assertion follows from Proposition \ref{prop:FN}. 
By induction on $d$ assume that we have a smooth $C_{d-1}$ ($d\geq
2$) satisfying (1)--(5).
We verify $C_{d-1}\cup l$ satisfies the following (1)'--(5)',
which are suitable modifications of (1)--(5):\\
(1)' $C_{d-1}\cup l$ intersects $B_{\varphi}$ simply by (1) for $C_{d-1}$
and generality of $l$.\\
(2)'
$M_{d-1}\cup M(l)$ intersects $Q_2$ simply by (2) for $C_{d-1}$
and generality of $l$.\\
(3)'
$M_{d-1} \cup M(l)$ is not irreducible but is of degree $d$
and has only simple nodes by (4) for $C_{d-1}$.\\
(4)'
$M_{d-1}\cup M(l)\cup M(m)$ has only simple nodes as its singularities
for a general line $m$ intersecting $C_{d-1}$.

Indeed, since $m$ is also general, 
$M_{d-1} \cup M(m)$ has only simple nodes by (4) for $C_{d-1}$.
Thus we have only to prove that
$M_{d-1}\cap M(l)\cap M(m)=\emptyset$, namely,
there is no secant line of $C_{d-1}$ intersecting both $l$ and $m$.
Fix a general $l$ and move $m$.
If there are secant lines $r_m$ of $C_{d-1}$ intersecting both $l$ and $m$
for general $m$'s,
then $r_m$ moves whence we have $M(l)\subset M_{d-1}$, a contradiction.
\\
(5)'
For a bi-secant line $\beta$ of $C_{d-1}\cup l$ except the lines
through $C_{d-1}\cap l$,
the curve $M_{d-1}\cup M(l)\cup M(\beta)$ has only simple nodes as its singularities.

Indeed, if $\beta$ is a bi-secant line of $C_{d-1}$, then the assertion follows from
(5) for $C_{d-1}$ by a similar way to the proof of (4)'.
Suppose that $\beta$ is a uni-secant line of $C_{d-1}$ intersecting $l$.
We have only to prove that
there is no secant line of $C$ intersecting both $l$ and $\beta$.
If there is such a line $r$, then
$l$, $\beta$ and $r$ pass through one point.
This does not occur for general $l$ and $\beta$ by
Proposition \ref{prop:Cd1} (5).

Thus, by a deformation theoretic argument, 
we see that $C_d$ satisfies (1)--(5). 
\end{proof}
\vspace{0.15in}

\subsection{On irreducibility of families of rational curves on $B$}~

We discuss about irreducibility of 
the Hilbert scheme of smooth rational curves on $B$ of
a fixed degree though we do not need it fully.
%See \cite{TZ} for further investigation.

%%%%%%%%UNICITA~%%%%%%%%%%%%%%%%%%%%%%%%%%%%%%%%%%%%%%%

For a smooth projective variety $X$ in some projective space,
let $\sH^{0}_{d}(X)$ be the Hilbert scheme of smooth rational curves on $X$
of degree $d$.
By \cite{Per}, $\sH^{0}_{d}(G(a,b))$ is non-empty and 
irreducible,
where $G(a,b)$ is the Grassmannian parameterizing
$a$-dimensional subvector spaces in a fixed
$b$-dimensional vector space. 

Let $\sH^{0'}_{d}(X)$ be 
the open subset of $\sH^{0}_{d}(X)$
parameterizing smooth rational curves on $X$
of degree $d$ with linear hull of maximal dimension.

Let $\sH_d^{B}$ be the Hilbert scheme of general smooth rational curves
on $B$ of degree $d$ 
obtained inductively as in Proposition $\ref{prop:Cd0}$. 

We can show inductively that
$\sH_d^{B}\subset {\sH}^{0'}_{d}(B)$, thus
we can ask the following:

\begin{que}
\label{prob:irred}
$\overline{\sH}_d^{B}=\overline{\sH}^{0'}_{d}(B)$ ?
(here we take the closure in the Hilbert scheme.)
Are they irreducible ?
\end{que}

We have a partial answer to this question as follows:

\begin{prop}
\label{prop:Cdbis} 
$\sH_d^{B}$ with any $d$ and
$\sH^{0'}_{d}(B)$ with $d\leq 6$ are irreducible.
$\overline{\sH}^{0'}_{d}(B)=\overline{\sH}_d^{B}$ for $d\leq 6$. 
\end{prop}

\begin{proof} 
%and, in particular, so is $\sH^{0}_{d}(\mP^{N})$. 
%Consider $G(2,5)\subset \mP^{9}$. 
%The claim is that $\sH^{0}_{d}(\mP^{6}\cap G(2,5))$ is
%irreducible, where $\mP^{6}\subset \mP^9$ is transversal to $G(2,5)$.
The claim is true for $d=1$ since $\overline{\sH}^{0'}_{1}(B)=\overline{\sH}^B_1\simeq \mP^2$.

First we prove 
$\sH_d^{B}$ is irreducible for any $d$.
By induction let us
assume that $\sH_{d-1}^B$ is irreducible. Let 
$[C^{0}_{d-1}]\in \sH_{d-1}^B$ be a generic element.
The family of lines $[l]\in\sH_{1}^{B}$ which intersect 
a general element of $\sH_{d-1}^B$ is irreducible 
by Proposition \ref{prop:Cd} (3). This implies that the family 
$\sH_{d-1,1}^B$ of reducible curves $C^{0}_{d}=C^{0}_{d-1}\cup l$
such that 
$[C^{0}_{d-1}]\in \sH_{d-1}^B$, $[l]\in\sH_{1}^{B}$ and
${\rm{length}}\, C^{0}_{d-1}\cap l=1$ is irreducible. 
As in the proof of Proposition \ref{prop:Cd0},
the Hilbert scheme is smooth at the point $[C^{0}_{d}]$.
Thus $\sH_d^{B}$ is irreducible.

Second we prove $\sH^{0'}_{d}(B)$ with $d\leq 6$ is irreducible.
Let $\sB$ be the irreducible family of del Pezzo 3-folds
$B=G(2,5)\cap\mP^{6}$, where 
$\mP^{6}\subset \mP^9$ is transversal to $G(2,5)$. 
Let 
$$
J=\{([C^{0}_{d}],[B])\in \sH^{0'}_d(G(2,5)) \times\sB \mid 
C^{0}_{d}\subset B\}.
$$
\noindent
%The claim is equivalent to show that 
%a general fiber $J\to \sB$ is irreducible. 
If $d\leq 6$, then it is known that
a general 
smooth rational curve of degree $d$ on $G(2,5)$ is a normal rational curve,
and is contained in a smooth $3$-dimensional linear section of $G(2,5)$,
namely, a smooth quintic del Pezzo $3$-fold.
Indeed, we can construct such a rational curve with $d\leq 5$ explicitly
on a smooth quintic del Pezzo surface, which is contained in
a smooth quintic del Pezzo $3$-fold.
For $d=6$, $C_{\varphi}$ as in the subsection \ref{subsection:FN}
is an example of such a rational curve $C_6$ on
a smooth quintic del Pezzo $3$-fold.

Thus a general fiber $J\to \sB$ is equal to
$\sH^{0'}_d(B)$ and is non-empty.
Moreover,
any fiber of $J\to \sH^{0'}_d(G(2,5))$ is isomorphic to
$G(\mP^d,\mP^6)$.
Since $\sH^0_d(G(2,5))$ is irreducible and 
$\sH^{0'}_d(G(2,5))$ is an open subset of $\sH^0_d(G(2,5))$,
it holds $J$ is irreducible.
By the argument of \cite[Proof of Theorem 3.1 p.17]{MT},
we have only to show that there is 
one particular component 
of a general fiber $J\to \sB$ invariant under monodromy.
Actually, this is nothing but $\sH^B_d$.
\end{proof}

\begin{cor}
\label{cor:hull}
Let $C_d$ be a general smooth rational curve 
constructed as in Proposition $\ref{prop:Cd0}$.
If $d=5$, 
then $C_5$ is a normal rational curve and 
is contained in a unique hyperplane section $S$, 
which is smooth.
If $d\geq 6$, then
$C_d$ is not contained in a hyperplane section.
\end{cor}

\section{Various projections of $B$}

%To investigate the geometry of $B$,
%we often project $B$ from a linear subspace.

\subsection{Projection of $B$ from a line or a conic}

%In this subsection,
%we describe various projection
%of $B$.

\begin{prop}
\label{prop:proj1}
\begin{enumerate}[$(1)$]
\item
Let $l$ be a line on $B$. Then the 
projection of $B$ from $l$ is decomposed as follows$:$
\begin{equation}
\label{eq:line}
\xymatrix{ & B_l \ar[dl]_{\pi_{1 l}} \ar[dr]^{\pi_{2 l}}\\
 B &  & Q, }
\end{equation}
where $\pi_{1l}$ is the blow-up along $l$ and 
$B \dashrightarrow Q$ is the projection from $l$ and 
$\pi_{2l}$ contracts onto a rational normal curve of degree $3$ the
strict transform of the loci swept by the lines of $B$
touching $l$. Moreover 
\begin{equation}
\label{eq:line2}
-K_{B_l}=H+H_Q,
\end{equation}
where $H$ and $H_Q$ are the pull backs of general hyperplane sections of 
$B$ and $Q$ respectively.
We denote by $E_l$ the $\pi_{1l}$-exceptional divisor.
\item
Let $q$ be a smooth conic on $B$. Then the 
projection of $B$ from $q$ behaves as follows$:$
\begin{equation}
\label{eq:conic}
\xymatrix{ & {B_{q}} \ar[dl]_{\pi_{1 q}} \ar[dr]^{\pi_{2 q}}\\
 B &  & \mP^{3} , }
\end{equation}
\noindent 
where $\pi_{1 q}$ is the blow-up of $B$ along $q$ and 
$\pi_{2 q}\colon {B_{q}}\rightarrow \mP^{3}$ is the
divisorial contraction of the strict transform $T_{q}$ of the loci
swept by
the lines touching $q$. 
%In particular if $E_{q}=¥\pi_{1 q}^{-1}(q)$
%then $T_{q}\in\mid 2H-3E_{q}\mid$ and 
Moreover
\begin{equation}
\label{eq:conic2}
-K_{{B_{q}}}=H+H_{\mP},
\end{equation}
where $H$ and $H_{\mP}$ are the pull backs of general hyperplane sections of 
$B$ and $\mP^3$ respectively.
\end{enumerate}
\end{prop}

\begin{proof} These results are more or less well-known. 
For (1), refer \cite{Fu2}, and 
for (2) (and (1)), refer \cite{MM1},
No. 22 for (2) (No. 26 for (1)).
See also \cite{MM3}, p.533 (7.7) for a discussion.
\end{proof}

We give several applications of the projection of $B$ from a line or a conic.

Let $C:=C_d$ be a general rational curve 
of degree $d$
constructed as in Proposition $\ref{prop:Cd0}$, and
${l}_1$ and ${l}_2$ two general secant lines of $C$
such that ${{l}_1}\cap {{l}_2}=\emptyset$.
We need to count the number of
multi-secant conics of $C$
intersecting $l_1$ and $l_2$
in the proof of Theorem \ref{thm:H_2}.

\begin{lem}
\label{lem:twolines}
Assume that $d\geq 3$.
Let $B\dashrightarrow Q \dashrightarrow \mP^2$
be the successive linear projections 
from ${l}_1$ and then
the strict transform of ${l}_2$ on $Q$.
Let ${l}$ be another general secant line of $C$,
and $C'$ and ${l}'\subset\mP^{2}$  be the images of $C$
and ${l}$ respectively. 
Then $C\cup {l}\dashrightarrow C'\cup {l}'$
is generically one to one and $\deg C'\cup l'=d-1$. 
Moreover,
$C'\cup {l}'$ has only simple nodes as its singularities.
In particular $($since $\deg C'=d-2$ and $C'$ is rational$)$
$C'$ has $\frac{(d-3)(d-4)}{2}$ simple nodes, equivalently,
there exist 
$\frac{(d-3)(d-4)}{2}$ bi-secant conics of $C$ intersecting both
${l}_1$ and ${l}_2$.
\end{lem}

\begin{rem}
The line $l$ is needed for the inductive proof as below.
\end{rem}

\begin{proof}
We show the assertion using the inductive construction of $C=C_d$.
The assertion follows for $d=3$ directly.
Consider a smoothing from $C_{d-1}\cup {m}$ to $C_d$.
Let ${m}_1$ and ${m}_2$ 
two general secant lines of $C_{d-1}$
such that ${{m}_1}\cap {{m}_2}=\emptyset$.
Let $B\dashrightarrow Q \dashrightarrow \mP^2$
be the successive linear projections 
from ${m}_1$ and then from
the strict transform of ${m}_2$ on $Q$.
Let ${r}$ be another general secant line of $C_{d-1}$,
and $C'_{d-1}, {m}'$ and ${r}' \subset \mP^2$
be the images of $C_{d-1}$, ${m}$ and ${r}$ respectively. 
Then we have only to show that
$C_{d-1}\cup {m}\cup {r} 
\dashrightarrow C'_{d-1}\cup {m}'\cup {r}'$ 
is generically one to one,
$\deg C'_{d-1}\cup {m}'\cup {r}'=d-1$ and 
$C'_{d-1}\cup {m}'\cup {r}'$ has 
only simple nodes as its singularities
assuming $C_{d-1}\cup {r}\dashrightarrow C'_{d-1}\cup {r}'$ 
is generically one to one,
$\deg C'_{d-1}\cup {r}'=d-2$ 
and $C'_{d-1}\cup {r}'$ has only simple nodes as its singularities.

Since ${m}$ is also general,
$C_{d-1}\cup {m}\dashrightarrow C'_{d-1}\cup {m}'$ is generically one to one,
$\deg C'_{d-1}\cup {m}'=d-2$
 and 
$C'_{d-1}\cup {m}'$ has only simple nodes as its singularities.
Thus $C_{d-1}\cup {m}\cup {r} 
\dashrightarrow C'_{d-1}\cup {m}'\cup {r}'$ 
is generically one to one and
$\deg C'_{d-1}\cup {m}'\cup {r}'=d-1$.
To show $C'_{d-1}\cup {m}'\cup {r}'$ has 
only simple nodes as its singularities,
it suffices to prove that
there are no secant conics of $C_{d-1}$
intersecting all the
${m}_1$, ${m}_2$, ${m}$ and ${r}$.
This follows from the fact that a secant conic ${q}$ of $C_{d-1}$
intersects finitely many secant lines of $C_{d-1}$
by $M({q})\not \subset M(C_{d-1})$.

The last statement follows from that,
by generality of $l_1$ and $l_2$,
any multi-secant conic of $C$ intersecting $l_1$ and $l_2$
is bi-secant.
\end{proof}

The following is a variant of Lemma \ref{lem:twolines},
which is also need in the proof of Theorem \ref{thm:H_2}.

\begin{lem}
\label{lem:twolines'}
Assume that $d\geq 4$.
Let ${l}_0$ be a general uni-secant line of $C$.
Let $B\dashrightarrow Q \dashrightarrow \mP^2$
be the successive linear projections 
from ${l}_0$ and then
the strict transform of a bi-secant line $\beta_i$ on $Q$.
Let ${l}$ be another general uni-secant line of $C$,
and $C'$ and ${l}'\subset\mP^{2}$  be the images of $C$
and ${l}$ respectively. 
Then $C\cup {l}\dashrightarrow C'\cup {l}'$ 
is generically one to one, $\deg C'\cup l'=d-2$, and 
$C'\cup {l}'$ has only simple nodes as its singularities. 
In particular $($since $\deg C'=d-3$ and $C'$ is rational$)$
$C'$ has $\frac{(d-4)(d-5)}{2}$ simple nodes,
equivalently,
there exist $\frac{(d-4)(d-5)}{2}$ bi-secant conics of $C$ intersecting
$\beta_i$ and ${l}_0$ except conics containing $\beta_i$.
\end{lem}

\begin{proof}
Similarly to the previous lemma, 
we show the assertion using the inductive construction of $C=C_d$.
The assertion follows for $d=4$ directly.
Consider a smoothing from $C_{d-1}\cup {m}$ to $C_d$.
Let ${m}_0$ be a general uni-secant line of $C_{d-1}$,
and $\beta$ a bi-secant line of $C_{d-1}\cup {m}$
different from any of the remaining two lines through $C_{d-1}\cap {m}$. 
Let $B\dashrightarrow Q \dashrightarrow \mP^2$
be the successive linear projections 
from ${m}_0$ and then
the strict transform of $\beta$ on $Q$.
Let ${r}$ be another general uni-secant line of $C_{d-1}$,
and $C'_{d-1}, {m}'$ and ${r}' \subset \mP^2$
be the images of $C_{d-1}$, ${m}$ and ${r}$ respectively. 

First we suppose that $\beta$ is a bi-secant line of $C_{d-1}$.
Then we have only to show that
$C_{d-1}\cup {m}\cup {r} 
\dashrightarrow C'_{d-1}\cup {m}'\cup {r}'$ 
is generically one to one, $\deg C'_{d-1}\cup {m}'\cup {r}'=d-2$,
and $C'_{d-1}\cup {m}'\cup {r}'$ has 
only simple nodes as its singularities
assuming $C_{d-1}\cup {r}\dashrightarrow C'_{d-1}\cup {r}'$ is birational and 
$C'_{d-1}\cup {r}'$ has only simple nodes as its singularities.
The proof is the same as that of Lemma \ref{lem:twolines}, 
so we omit it.

Next suppose that $\beta$ is a uni-secant line of $C_{d-1}$ 
intersecting ${m}$ outside $C_{d-1}\cap {m}$.
Note that, by the projection $B\dashrightarrow \mP^2$,
${m}$ is contracted to a point.
Moreover,
$\beta$ is a general uni-secant line since so is ${m}$.
Thus, by Lemma \ref{lem:twolines}, 
$C_{d-1}\cup {m}\cup {r} 
\dashrightarrow C'_{d-1}\cup {r}'$ 
is generically one to one, $\deg C'_{d-1}\cup {r}'=d-2$,
and 
$C'_{d-1}\cup {r}'$ has 
only simple nodes as its singularities.
\end{proof}

Let $f\colon A\to B$ be the blow-up of $B$ along a general
smooth rational curve $C_d$. 
The following lemma can be regarded as the assertion of generality of $C_d$.
We need this in the subsection \ref{subsection:finite}.

\begin{lem}\label{sullebisecanti}
Let $\beta'_i\subset A$ be the strict transform of a bi-secant line $\beta_i$
of $C_d$.
It holds:
    $$
    \sN_{\beta'_{i}/A}=\sO_{\beta'_{i}}(-1)\oplus\sO_{\beta'_{i}}(-1).
    $$
\end{lem}

\begin{proof}
We prove the assertion by using the inductive construction of $C_d$.
The assertion is clear for $d=1$ since $C_1$ has no bi-secant line.

Suppose the assertion holds for $C_{d-1}$. 
Choose a general uni-secant line $l\subset B$ of $C_{d-1}$.
Let ${m}_1,\dots, {m}_{d-2}$ be the lines on $B$ 
intersecting both $C_{d-1}$ and ${l}$ outside 
$C_{d-1}\cap {l}$. By generality of $C_{d-1}$ we can assume that 
${m}_1,\dots, {m}_{d-2}$ are unisecant of $C_{d-1}$.

Let $A'\to B$ be the blow-up along $C_{d-1}\cup {l}$.
Note that the smoothing $C_{d-1}\cup {l}$ to $C_d$ induces 
that of $A'$ to $A$.
Let $\widetilde{m}_i$ be
the strict transform of ${m}_i$ on $A'$.
By the smoothing construction of $C_d$ from $C_{d-1}\cup {l}$
and the assumption on induction, 
we have only to prove $\sN_{\widetilde{m}_i/A'}=
\sO_{\mP^1}(-1)\oplus\sO_{\mP^1}(-1)$.
Let $A'_1\to B$ be the blow-up along ${l}$ and 
$A'_2\to A'_1$ the blow-up along the strict transform of $C_{d-1}$.
Denote by $m'_i$ and $m''_i$ the strict transform of ${m}_i$ on 
$A'_1$ and $A'_2$ respectively.
Then $\sN_{\widetilde{m}_i/ A'}=\sN_{m''_i/ A'_2}$.
We consider the projection of $B$ from the line $l$ as in
Proposition \ref{prop:proj1} (2). 
Since $m'_i$ is a fiber of $A'_1\to Q$, 
we have $\sN_{m'_i/A'_1}=\sO_{\mP^1}\oplus\sO_{\mP^1}(-1)$.
Let $F$ be the exceptional divisor of $A'_1\to Q$ and $F'$ the strict transform
of $F$ on $A'_2$. We may suppose $F$ and $C'_{d-1}$ intersect transversely,
thus $F'\to F$ is the blow-up at $d-2$ points $m'_i\cap C'_{d-1}$ 
($i=1,\dots, d-2$). Thus $F'\cdot m''_i=-1$ and $\sN_{m''_i/F'}=
\sO_{\mP^1}(-1)$, and this implies the assertion.
\end{proof}

\subsection{Double projection of $B$ from a point}
\label{subsection:doubleproj}
\begin{defn}
Let $b$ be a point of $B$.
We call the rational map from $B$
defined by the linear system of hyperplane sections
singular at $b$ {\em{the double projection from $b$}}.
\end{defn}

\begin{prop}
\label{prop:proj2}
Let $b$ be a point of $B$. Then
\begin{enumerate}[$(1)$]
\item
 
the target of the double projection from $b$ is $\mP^2$,
and the double projection from $b$ and 
the projection $B\dashrightarrow \overline{B}_b$ from $b$ 
fit into the following diagram:
\begin{equation}
\label{eq:pt1}
\xymatrix{ & B_b \ar[dl]_{\pi_{1 b}} \ar[dr] & & 
 B'_b \ar[dl] \ar[dr]^{\pi_{2 b}} & \\
 B &  & \overline{B}_b & & \mP^2, }
\end{equation}
where $\pi_{1 b}$ is the blow-up of $B$ at $b$,
$B_b\dashrightarrow B'_b$
is the flop of the strict transforms of lines through $b$, 
and $\pi_{2 b}\colon B'_b\rightarrow \mP^{2}$ 
is a $(\text{unique})$ $\mP^1$-bundle structure.
Moreover, $\overline{B}_b\dashrightarrow \mP^2$ is
the projection from the plane
which is the image of $\pi_{1b}$-exceptional divisor.

We denote by
$E_b$ the $\pi_{1b}$-exceptional divisor and
by $E'_b$ the strict transform of $E_b$ on $B'_b$,
%Moreover we have the following descriptions$:$
\item 
\label{eq:pt2}
\[
L=H-2E'_b\ \text{and}\ -K_{B'_b}=H+L,
\]
where 
$H$ is the strict transform 
of a general hyperplane section of $B$,
and $L$ is the pull back of a line on $\mP^2$,
\item
Case $(\mathrm{a})$\\
If $b\not \in B_{\varphi}$,
then the strict transforms $l'_i$ of three lines $l_i$ through $b$ on $B_b$
have the normal bundle $\sO_{\mP^1}(-1)\oplus \sO_{\mP^1}(-1)$.
The flop $B_b\dashrightarrow B'_b$ is the Atiyah flop.
In particular, $E'_b\to E_b$ is the blow-up at
the three points $E_b\cap l'_i$.\\
Case $(\mathrm{b})$\\
If $b \in B_{\varphi}\setminus C_{\varphi}$, then 
$E_b\dashrightarrow E'_b$ can be described as follows$:$
let $l$ and $m$ be two lines through $b$, where
$l$ is special, and $m$ is not special.
Let $l'$ and $m'$ be the strict transforms of $l$ and $m$ on $B_b$.
First blow up $E_b$ at two points $t_1:=E_b\cap l'$ and $t_2:=E_b\cap m'$
and then blow up at a point $t_3$ on the exceptional curve $e$ over $t_1$. 
Finally, contract the
strict transform of $e$ to a point. Then we obtain $E'_b$
$($this is a degeneration of the case $(\mathrm{a})$$)$.\\
Case $(\mathrm{c})$\\
See $\textup{\cite{FuNa2}}$ in case of $b \in C_{\varphi}$, and 
%\item
%$b \in B_{\varphi}\setminus C_{\varphi}$.
%Let $l$ and $m$ be two lines through $b$, where
%$l$ is special, and $m$ is not special.
%The flop of the strict transform of $m$ is the Atiyah flop as in case 
%$(\mathrm{a})$. Let $m'$ be the strict transform of $m$ on $B_b$.
%Then the flop of $m'$ is constructed as follows:
%let $T'\subset B_b$ be the strict transform of $T_m$, 
%the locus of lines intersecting $m$. Let $B_1\to B_b$ be the blow-up along $m'$%,and
%$T_1\subset B_1$ the strict transform of $T_m$. Then
%$T_1$ is the blow-up of the normalization of $T_m$ at two points
%$t_1$ and $t_2$, where $t_1$ and $t_2$ are contained in the
%inverse image of the non-normal locus, $t_1$ is on the negative section,
%and $t_2$ is on the fiber. 
%Let $m_1\subset T_1$ be the strict transform of the negative section
%and $\gamma\subset T_1$ the strict transform of the fiber containing $t_2$.
%Then $N_{m_1/B_1}$
%\item
%$b \in C_{\varphi}$.
%See \cite{FuNa2},
%and
%\end{enumerate}
\item
a fiber of $\pi_{2 b}$ not contained in $E'_b$
is the strict transform of a conic through $b$, 
%(which is disjoint from the strict transforms of lines through $b$), 
or
the strict transform of a line $\not \ni b$ intersecting a line through $b$.

The description of the fibers of $\pi_{2 b}$ contained in $E'_b$ is 
as follows$:$\\
Case $(\mathrm{a})$\\
If $b\not \in B_{\varphi}$, then
$\pi_{2 b|E'_b}\colon E'_b \to \mP^2$ is the blow-down
of the strict transforms of three lines
connecting two of $E_b\cap l'_i$, namely,
$E_b\dashrightarrow \mP^2$ is the Cremona transformation.\\
Case $(\mathrm{b})$\\
Assume that $b \in B_{\varphi}\setminus C_{\varphi}$. 
Then $\pi_{2 b|E'_b}\colon E'_b \to \mP^2$ is the blow-down
of the strict transforms of two lines, one is the line 
connecting $t_1$ and $t_2$,
the other is the line whose strict transform passes through $t_3$. 
$E_b\dashrightarrow \mP^2$ is a degenerate Cremona transformation.\\
Case $(\mathrm{c})$\\ 
See $\textup{\cite{FuNa2}}$ in case of $b \in C_{\varphi}$. 
\end{enumerate}
\end{prop}

\begin{proof}
This is a standard result 
in the birational geometry of Fano $3$-folds but 
is less known than Proposition \ref{prop:proj1}. 
We have only found the paper \cite{FuNa2},
in which they deal with the most difficult case (c).
Here we sketch 
the construction of the flop in the middle case (b) to intend the reader
to get a feeling of birational maps from $B$.

Let $b$ be a point of $B_{\varphi}\setminus C_{\varphi}$.
We use the notation of the statement of (3). 
The flop of $m'$ is the Atiyah flop. We describe the flop of $l'$.
By $\sN_{l/B}\simeq \sO_{\mP^1}(1)\oplus \sO_{\mP^1}(-1)$,
it holds that $\sN_{l'/B_b}\simeq \sO_{\mP^1}\oplus \sO_{\mP^1}(-2)$.
Hence the flop of $l'$ is a special case of Reid's one \cite[Part II]{Pagoda}.
We show that the width is two in Reid's sense.  
Let $T_1$ be the normalization of $T_l$. By Proposition \ref{prop:FN} (5), 
$T_1\simeq \mF_3$ and
the inverse image of the singular locus of $T_l$ is the union of
the negative section $C_0$ and a fiber $r$.
Let $\mu\colon \widetilde{B}_b\to B_b$ be the blow-up along $l'$ and 
$F$ the exceptional divisor. Let $T_2$ be the strict transform of $T_l$ on 
$\widetilde{B}_b$.
Then $T_2$ is the blow-up of $T_1$ at two points $s_1\in C_0$ and $s_2\in r$. 
Denote by $C'_0$ and $r'$ the strict transforms of $C_0$ and $r$.
We prove that $\sN_{r'/\widetilde{B}_b}\simeq \sO_{\mP^1}(-1)^{\oplus 2}$.
Note that $F\cap T_2=C'_0\cup r'$. The curves $C'_0$ and $r'$ are two sections
on $F$. Let $T'_1$ be the image of $T_2$ on $B_b$.
By $\sN_{l'/B_b}\simeq \sO_{\mP^1}\oplus \sO_{\mP^1}(-2)$
and $T_2=\mu^* T'_1-2F$,
it holds $F\simeq \mF_2$, and $T_{2|F}\sim 2G_0+3\gamma$,
where $G_0$ is the negative section of $F$ and 
$\gamma$ is a fiber of $F\to l'$.
Note that $F\cdot C'_0=(F_{|T_2}\cdot C'_0)_{T_2}=-3$ and 
$F\cdot r'=(F_{|T_2}\cdot r')_{T_2}=0$,
and 
$F\cdot G_0=0$ and $F\cdot (G_0+3\gamma)=-3$.
Thus we have $C'_0\sim G_0+3\gamma$ and $r'=G_0$ on $F$. 
Now we see that $-K_{\widetilde{B}_b}\cdot r'=
(\mu^*(-K_{B_b})-F)\cdot r'=0$.
Therefore, by $(r')^2=-1$ on $T_2$, it holds that
$\sN_{r'/\widetilde{B}_b}\simeq \sO_{\mP^1}(-1)^{\oplus 2}$.

It is easy to see that we can flop $r'$.
Let $\widetilde{B}_b\dashrightarrow \widetilde{B}'_b$ be the
flop of $r'$ (now we consider locally around $r'$).
Let $F'$ be the strict transform of $F$ on $\widetilde{B}'_b$.
By \cite{Pagoda}, $F'\simeq F$ and 
there is a blow-down $\widetilde{B}'_b\to \widetilde{B}''_b$ of $F'$ 
such that $\widetilde{B}''_b$ is smooth.
$\widetilde{B}_b\dashrightarrow \widetilde{B}''_b$ is the flop of $l'$.

By this description of the flop, 
we can easily obtain (3).
\end{proof}

As a first application of the above operations, we have the following result,
which we often use: 

\begin{cor}
\label{cor:uniqueconic}
Let $b_1$ and $b_2$ be two $($possibly infinitely near$)$ points on
$B$
such that there exists no line on $B$ through them.
Then there exists a unique conic on $B$ through $b_1$ and $b_2$.
\end{cor}

\begin{proof}
We project $B$ from $b_1$ as in (\ref{eq:pt1}).
Then the assertion follows by the description of fibers
of ${\pi_{2 b_{1}}}$ as in Proposition \ref{prop:proj2} (4). 
\end{proof}

\begin{nota}
\label{nota:Cb}
Consider the double projection from $b$, see proposition \ref{prop:proj2}.
Throughout the paper,
we denote by $C'_b$, $C''_b$ and $C_b$ 
the strict transforms of $C:=C_d$ on $B_b$, $B'_b$ and $\mP^2$
respectively. 
\end{nota}

The following result is one of the key results
for the proof of the main result.
Its importance and difficulty
lies in the actual fact that it holds not only for a general $b\in B$
but also for
every $b\in B$.

\begin{prop}
\label{prop:Cd2}
Let $C_d$ be a general smooth rational curve of 
degree $d$ on $B$ constructed as in Proposition $\ref{prop:Cd0}$.
Assume that $d\geq 5$.
Then,
for any point $b\in B$,
the restriction of $\pi_b$ to $C_d$ is birational.
\end{prop}
\begin{proof}
We prove the assertion by induction
based on the construction of $C_d$ from 
$C_{d-1}\cup l$, where
$l$ is a general uni-secant line of $C_{d-1}$ on $B$.

First we prove the assertion for $d=5$.
Assume by contradiction that 
$\pi_{b|C_5}$ is not birational for a point $b$.
Then, since $C\dashrightarrow C_b$ is a composite of linear projections,
$C_b$ is a line or conic in $\mP^2$.
Let $S$ be the pull-back of $C_b$ by $\pi_{2b}$.
If $C_b$ is a line, then $C_5$ is contained in a
singular hyperplane section, which is the strict transform of
$S$ on $B$ 
(recall that $B\dashrightarrow \mP^2$ is the double projection from $b$). 
This contradicts Corollary \ref{cor:hull}.
Assume that $C_b$ is a conic.
The only possibility is that
$L\cdot C''_b=4$ and $C''_b\to C_b$ is a double cover
since $L\cdot C''_b =\deg C_b \cdot \deg (C''_b\to C_b)\leq 5$. 
Since the flop does not change the intersection numbers 
between the canonical divisor and curves,
we have $-K_{B'_b}\cdot C''_b=-K_{B_b}\cdot C'_b$.
If $b\in C$, then we have $-K_{B'_b}\cdot C''_b=8$.
Thus, by Proposition \ref{prop:proj2} (2) and
$L\cdot C''_b=4$, it holds $H\cdot C''_b=4$.
By $L=H-2E'_b$, this shows that $E'_b\cdot C''_b=0$.
This is, however, a contradiction since
$E'_b\cap C''_b\not =\emptyset$.
Thus $b\not \in C$, and, by Proposition \ref{prop:proj2} (2),
it holds $H\cdot C''_b=6$.
By $L=H-2E'_b$, we have $E'_b\cdot C''_b=1$.
We compute ${E'_b}^2 S$.
Note that $-K_{B'_b}=2H-2E'_b=2(L+2E''_b)-2E''_b=2(L+E''_b)$.
We have
\[
{E'_b}^2 L=\frac 14 (-K_{B'_b}-2L)^2 L=
\frac 14 (-K_L-L_{|L})^2=1.
\]
Thus we have ${E'_b}^2 S=2{E'_b}^2 L=2$. 
The surface $S$ is a Segre-del Pezzo scroll.
Let $C_0$ is the negative section of $S$ and $l$ is a fiber of $S\to C_b$
and set $e:=-C^2_0$.
We can write 
$E'_{b|S}\sim C_0+pl$ and $C''_b\sim 2C_0+ql$ $(p,q\geq 0)$,
By $E'_b\cdot C''_b=1$ and ${E'_b}^2 S=2$,
we have $q+2p-2e=1$ and $2p-e=2$.
Thus $e=2p-2$ and $q=2p-3$.
Since $C''_b$ is irreducible,
$q\geq 2e$, whence $2p-3\geq 2(2p-2)$, i.e., $p=0$ and $q=-3$, a contradiction.

Assume that $d\geq 6$.
Let $\sC\to \Delta$ be the one-parameter smoothing of $C_{d-1}\cup l$
such that $\sC$ is smooth
(as we saw in the proof of Proposition \ref{prop:Cd0}, this is possible).
We consider the trivial family of the double projections
$B\times \Delta\dashrightarrow \mP^2\times \Delta$
from $b\times \Delta$.
Denote by $\sC'_b, \sC''_b$ and $\sC_b$ the strict transforms of
$\sC$ on $B_b\times \Delta$, $B'_b\times \Delta$
and $\mP^2\times \Delta$ respectively.
We also denote by $C'_{d-1,b}$, $C''_{d-1,b}$, and $C_{d-1,b}$
the strict transforms of
$C_{d-1}$ on $B_b$, $B'_b$ and $\mP^2$ respectively.
To prove the proposition,
it suffices to
show that, for any $b$, there exists at least one point on 
$C_{d-1,b}$ over which $\sC\dashrightarrow \sC_b$ is
isomorphic. First, admiting this claim, 
we finish the proof of the proposition.
Indeed, set 
\[
\sN:=
\{(b, t)\in B\times \Delta\mid \sC\dashrightarrow \sC_b 
\ \text{is not isomorphic over any point of } \ \sC_{b,t} \}
\]
and let $\Delta'\subset \Delta$ be the image of $\sN$ 
by the projection to $\Delta$. 
$\sN$ is a closed subset, and so is $\Delta'$ since $B\times \Delta\to \Delta$
is proper. 
Thus $\Delta'$ consists of finitely many points 
since the origin is not contained in $\Delta'$ by
admiting the above claim.
Therefore,
for a point $t\in \Delta$ sufficiently near the origin,
$\sC_t\dashrightarrow \sC_{t,b}$ is birational for any $b$,
which implies the proposition.

Now we show the above claim.
By induction, we may assume that $C_{d-1}\dashrightarrow C_{d-1,b}$
is birational for any $b$.
Note that $C_{d-1,b}$ is not a line since
otherwise $C_{d-1}$ is contained in 
a singular hyperplane section as we see above in the case of $C_5$,
a contradiction. 
We investigate the image of $l$ on $\mP^2$.
Recall the description of the fibers of $\pi_{2b}$ (Proposition \ref{prop:proj2} (4)).
If $b\not \in l$, then
the image of $l$ is a line or a point on $\mP^2$.
If $b\in l$, then the strict transform of $l$ on $B_b$ is a flopping
curve. Thus $\sC_b$ contains the image of the flopped curve,
which is a line.
We investigate the other possible irreducible components of 
the central fiber $\sC_{b,0}$ of $\sC_b\to \Delta$.  
If $b\not \in C_{d-1}\cup l$, then
the only possibility is that 
$\sC_{b,0}$ contains the image of a flopped curve,
which is a line on $\mP^2$.
%Thus 
%$\sC\dashrightarrow \sC_b$ is birational at a point of $C_{d-1}$.
Suppose $b\in C_{d-1}\cup l$.
Let $m'_b$ be the exceptional curve for $\sC'_b\to \sC$.
Since $\sC$ is a smooth surface,
$m'_b$ is a line on $E_b$.
The curve $\sC_{b,0}$ contains 
the strict transform $m_b$ of $m'_b$.
This is the only possibility of the other components of 
$\sC_{b,0}$.
Let $l'_b$ be the strict transform of $l$ on $B_b$ 
If $b\in l$, then
by the description of $E_b\dashrightarrow \mP^2$,
$m_b$ is a line since $l'_b$ is a flopping curve.
%Thus 
%$\sC\dashrightarrow \sC_b$ is birational at a point of $C_{d-1}$.
Suppose that $b\in C_{d-1}\setminus l$.
If $m'_b$ intersects a flopping curve,
$m_b$ is a line or a point.
%thus 
%$\sC\dashrightarrow \sC_b$ is birational at a point of $C_{d-1}$.
In the other case, $m_b$ is a conic.
If $b\not \in \cup_i \beta_i$, then
$\deg C_{d-1,b}=d-3$ by Proposition \ref{prop:proj2} (2).
By $d\geq 6$, $C_{d-1,b}$ is not a conic. 
Thus $C_{d-1,b}\not =m_b$. 
Assume $b\in \beta_i$.
Then 
$\deg C_{d-1,b}=d-4$. 
Thus, if $d\geq 7$, then 
$C_{d-1,b}\neq m_b$.
We show that
even if $d=6$, it holds $C_{d-1,b}\neq m_b$.
By Proposition \ref{prop:Cd1} (4),
the flop $B_b\dashrightarrow B'_b$ is of type (a) in Proposition \ref{prop:proj2} (3).
The strict transform $m''_b$ of $m'_b$ on $B''_b$
intersects the three fibers of $\pi_b$ contained in $E'_b$,
which are the strict transforms of three lines on $E_b$.
On the other hand, by $E'_b\cdot C''_{d-1,b}=2$,
the curve $C''_{d-1,b}$ intersects
at most two fibers of $\pi$ contained in $E'_b$.
Thus it holds $C_{d-1,b}\neq m_b$.

The above investigation shows that 
$\sC\dashrightarrow \sC_b$ is isomorphic over a point of $C_{d-1,b}$.
\end{proof}

We restate the proposition in terms of the relation 
between $C_d$ and multi-secant conics of $C_d$ on $B$ as follows:
\begin{cor}
\label{cor:finite}
Let $b$ be a point of $B$ not in any bi-secant line of $C_d$ on $B$.
If $d\geq 5$, then
there exist finitely many $k$-secant conics of $C_d$ on $B$ 
through $b$ with $k\geq 2$ if $b\not \in C_d$
$(\text{resp. with $k\geq 3$ if $b \in C_d$})$.
\end{cor}

\begin{proof}
For a point $b\in B$ outside bi-secant lines of $C_d$ on $B$,
there exist a finite number of 
singular multi-secant conics of $C_d$ through $b$
since the number of lines through $b$ is finite, and
the number of lines intersecting both a line through $b$ and $C_d$
is also finite by Proposition \ref{prop:Cd} (3).
Therefore
we have only to consider smooth multi-secant conics $q$ of $C_d$ through $b$.
By Proposition \ref{prop:proj2} (4),
the strict transform $q'$ of 
such a conic $q$ on $B'_b$ is a fiber of $\pi_{2b}$.
If $b\not \in C_d$,
then $q'$ intersects $C'_b$ twice or more counted with multiplicities, 
thus by Proposition \ref{prop:Cd2}, the finiteness of such a $q$ follows.
We can prove the assertion in case of $b\in C_d$ similarly, 
thus we omit the proof. 
\end{proof}

\begin{rem}
We refine this statement in Lemmas \ref{fuoripiatto} and
\ref{finitezzafinale}.
\end{rem}

\begin{lem}
\label{lem:Cb}
Let ${l}$ be a general uni-secant line of $C$ and
$l_b\subset \mP^2$ the image of $l$ by the double projection from a point $b$.
For a general point $b\not \in C$,
$\deg C_b=d$ and $C_b\cup l_b$ has only simple nodes.
Assume that $d\geq 3$.
For a general point $b$ of $C$,
$\deg C_b=d-2$ and $C_b\cup l_b$ has only simple nodes.
\end{lem}

\begin{proof}
The claims for $\deg C_b$ follows from Propositions 
\ref{prop:proj2} (2) and \ref{prop:Cd2}.
As for the singularity of $C_b\cup l_b$, the claim follows from
simple dimension count. For simplicity, we only prove that
for a general point $b\not \in C$, the curve $C_b$ has only simple nodes.
By Proposition \ref{prop:Cd1'}, we may assume that 
any multi-secant conic through $b$ is smooth, bi-secant and
intersects $C$ simply.
Let ${q}$ be a smooth bi-secant conic through $b$.
We may assume that $\sN_{{q}/B}\simeq \sO_{\mP^1}(1)^{\oplus 2}$.
Let $q'$ be the strict transform of ${q}$ on $B'_b$. 
Let $\widetilde{B}'\to B'_b$ be the blow-up along $q'$,
$E_{q'}$ the exceptional divisor and $\widetilde{C}''$ the strict transform of
$C''_b$. 
Note that $E_{q'}\simeq \mP^1\times \mP^1$
since $\sN_{q'/B'_b}\simeq \sO_{\mP^1}^{\oplus 2}$.
Then $C_b$ has simple nodes at the image of $q'$ if and only if
the two points in $E_{q'}\cap \widetilde{C}''$ does not belong to the same
ruling with the opposite direction to a fiber of $E_{q'}\to q'$. 
Let $\widetilde{B}_{{q}}\to B$ be the blow-up along ${q}$,
$E_{{q}}$ the exceptional divisor and $\widetilde{C}$ the
strict transform of $C$.
It is easy to see that 
a ruling of $E_{{q}}$ 
with the opposite direction to a fiber of $E_{{q}}\to {q}$
corresponds to that 
of $E_{q'}$ with the opposite direction to a fiber of $E_{q'}\to q'$.
Thus
$C_b$ has simple nodes at the image of $q'$ if and only if
the two points in $E_{{q}}\cap \widetilde{C}$ does not belong to the same
ruling with the opposite direction to a fiber of $E_{{q}}\to {q}$. 
We can show that this is the case for a general $b$ by simple dimension count.
\end{proof} 

\begin{cor}
\label{cor:n}
\begin{enumerate}[$(1)$]
\item
The number of multi-secant conics of $C$
through a general point of $B$ is $n:=\frac{(d-1)(d-2)}{2}$.
\item
The number of $k$-secant conics of $C$ with $k\geq 3$
through a general point of $C$ is $\frac{(d-3)(d-4)}{2}$.
\item
Let $l$ be a general uni-secant line of $C$.
Then the number of multi-secant conics of $C$
intersecting $l$
and passing through a general point of $C$
is $d-3$.
\end{enumerate}
\end{cor}

\begin{proof}
We only prove (1) since the other statement can be proved
similarly.

Let $b\not \in C$ be a general point of $B$.
Recall that, by
Corollary \ref{cor:finite},
there exist only finitely many multi-secant conics of
$C$ through $b$.
Moreover,
since $C_b$ is a nodal rational curve of degree $d$
by Lemma \ref{lem:Cb}, 
the number of its nodes is exactly $n$,
which is nothing but the number of 
multi-secant conics through $b$.
\end{proof}

As we saw in Corollary \ref{cor:n} (1),
a general point of $B$ gives $n$ multi-secant conics of $C$ 
through it.
Conversely, we ask whether mutually intersecting $n$ multi-secant conics of $C$
actually pass through one point or not.
The next lemma partially answer this question and
it is sufficient for our purpose in the proof of Theorem \ref{diretto}.
We remark that the case $d=5$ is treated in \cite[4.3]{doldual}.

%%%%%%%%%%%%%%%%%%%%%%%%%%%%%%%%%%%%%%%%%%%%%%%%%%%%%%%%%%%%%%%%%%%%%%%%%%%%%%%%
\begin{lem}
\label{lem:GeomA}
Let $q_1,\dots, q_n$ be mutually intersecting 
$n$ distinct multi-secant conics of $C$ 
such that
\begin{enumerate}[$(1)$]
\item
all ${q}_i$ are smooth, 
\item
no two of ${q}_i$ intersect at a point of $C\cup \cup_i \beta_i$, and
\item
if three of ${q}_i$ pass through a point $b$,  
then any other ${q}_i$ does not intersect a line through $b$
outside $b$.
\end{enumerate}
Then all $q_i$ pass through one point.
\end{lem}

\begin{rem}
The set of $n$ conics through a general point satisfies
the conditions of the lemma.
\end{rem}

\begin{proof}~

\noindent
{\bf{Step 1.}}
Let $b\in B$ be a point such that
five of ${q}_i$, 
say, ${q}_1,\dots, {q}_5$ pass through $b$.
%and no line through $b$ does not intersect any ${q}_j$ outside $b$.
Then all the ${q}_i$ pass through $b$.

By the double projection from $b$,
${q}_1,\dots,{q}_5$ are mapped to 
points $p_1,\dots,p_5$ on $\mP^2$.
%and  nor intersect any line th%rough $b$,
Suppose by contradiction 
that a smooth conic ${q}_j$ does not pass through $b$.
Let $q'_j$, $q''_j$ and $\widetilde{q}_j$ be the strict transforms
of ${q}_j$ on $B_b$, $B'_b$ and $\mP^2$, and set 
$S:=\pi_{2b}^*\widetilde{q}_j$. 
By the assumption (3), ${q}_j$ does not intersect
a line through $b$.
Thus $\widetilde{q}_j$ is a smooth conic 
through $p_1,\dots,p_5$. The conic $\widetilde{q}_j$ is unique
since a conic through five points is unique.
It holds that
$-K_{B'_b}\cdot q''_j=4$ and $S\cdot q''_j=4$, thus
$S\simeq \mF_2$ and $q''_j$ is the negative section.
This implies that $q_j$ is also unique.
%\footnote{Only the image of $q_j$ is unique...} 
By reordering, we may assume that $j=n$.
We have the configuration such that
all the conics pass through $b$ except $q_n$.
Denote by $p_i$ the image of $q_i$ ($i\neq n$).
Then $\widetilde{q}_n$ and $C_b$ intersect at $p_i$.
By $d\geq 6$, it holds $\deg C_b\geq 3$, thus
$\widetilde{q}_n\neq C_b$. 
By the assumption (2), $b\not \in C$.
Therefore $\widetilde{q}_n$ and $C_b$ intersect
at $n-1$ singular points of $C_b$.
Since $\deg C_b\leq d$, it holds
$2(n-1)\leq 2d$, a contradiction.\\
{\bf{Step 2.}}
If four conics ${q}_1, \dots, {q}_4$ 
pass through one point $b$, 
%and no line through $b$ does not intersect any $\overline{q}_j$ outside $b$,
then all the conics pass through $b$.

By contradiction and Step 1,
we may assume that all the conics except 
${q}_1,\dots,{q}_4$
do not pass through $b$.
Pick up two any conics, say, ${q}_5$ and ${q}_6$, 
not passing through $b$.
Considering the double projection from $b$ as in Step 1.
Denote by $\widetilde{q}_j$ $(j\geq 5)$ 
the image of ${q}_j$ on $\mP^2$. 
By the assumption (3), ${q}_5$ and ${q}_6$
do not intersect a line through $b$,
thus $\widetilde{q}_5$ and $\widetilde{q}_6$ are conics on $\mP^2$.
Therefore ${q}_5 \cap {q}_6$ lies 
on one of ${q}_1,\dots, {q}_4$
since otherwise 
$\widetilde{q}_5$ and $\widetilde{q}_6$ would intersect at five points
and this is a contradiction as in Step 1.
Thus any two conics intersect on ${q}_1,\dots, {q}_4$.
Let $p_i$ be the intersection 
${q}_i\cap {q}_5$ for $i=1,\dots, 4$.
Then ${q}_j$ ($j\geq 5$) pass through one of $p_i$.
Thus one of $p_i$, say, $p_1$, there
pass through at least $\lceil \frac{(n-5)}{4} \rceil$ conics.
By Step 1,
$\lceil \frac{(n-5)}{4}\rceil \leq 2$ 
(already ${q}_1$ and ${q}_5$ pass through $p_1$).
This implies $d=6$. We exclude this case in Step 3.
Note that if $d=6$, then the four conics ${q}_1$,
${q}_2$, ${q}_5$, and ${q}_6$
mutually intersect and the all the intersection points are different.
By reordering conics, we assume that 
${q}_i$ $(1\leq i \leq 4)$ satisfy this property in Step 3.\\ 
{\bf{Step 3.}} 
We complete the proof.

Assume by contradiction that
${q}_1,\dots,{q}_n$  do not pass through one point on $B$. 
If $d\geq 7$, then, by Steps 1 and 2,
\begin{equation}
\label{eq:three}
\text{
at most three of ${q}_i$'s pass through any intersection point.}
\end{equation}
Let $m$ be the number of conics in a maximal tree $T$ of ${q}_i$'s 
such that two conics in $T$ pass through any intersection point. 
Note that $T$ is connected since ${q}_i$'s mutually 
intersect. 
The number of the intersection points of ${q}_i$'s contained in 
$T$ is $\frac{m(m-1)}{2}$.

By the maximality of $T$,
a conic not belonging to $T$ passes through one of the
intersection points of conics in $T$.
By (\ref{eq:three}),
no two conics not belonging to $T$ 
pass through one of the intersection point of conics in $T$.
Hence it holds $\frac{m(m-1)}{2} +m \geq n$.
This implies that $m\geq d-2$ by $n=\frac{(d-1)(d-2)}{2}$.
By reordering, we assume that 
${q}_1,\dots, {q}_m$ belong to $T$.
If $d=6$, then we take ${q}_1,\dots, {q}_4$
as in the last part of Step 2.
Consider the projection $B \dashrightarrow \mP^3$ from the conic ${q}_1$.
Then ${q}_2,\dots,{q}_m$ are mapped 
to lines $l_2,\dots, l_m$ intersecting
mutually on $\mP^3$ and the intersection points are different.
Thus  $l_2,\dots, l_m$ span a plane, which in turn shows that
${q}_1,\dots, {q}_m$ span a hyperplane section $H$ on $B$.
Since $C$ intersects ${q}_i$ at two point or more, 
$C$ intersects $H$ at $2m$ points or more by the assumption (2).
But $2m\geq 2(d-2) >d$, $C$ must be contained in $H$,
a contradiction to Corollary \ref{cor:hull}. 
\end{proof}

%The results in \ref{subsection:CdB} are delicate but their proof is
%more or less based on standard parameter count.

\section{Lines and conics on $A$}
We fix a general $C:=C_{d}$ as in the subsection
\ref{subsection:CdB}.
Let $f\colon A\to B$ be the blow-up along $C$.
We start the study of the geometry of $A$. 
In the subsections \ref{subsection:H1} and \ref{subsection:H2},
we study the families of curves on $A$ of degree one or two with respect to
the anti-canonical sheaf of ${A}$ (we call them {\em{lines}} and {\em{conics}} on $A$ respectively).
The curve $\sH_1$ parameterizing lines on $A$ and 
the surface $\sH_2$ parameterizing conics on $A$ are two of the main characters
in this paper. See Corollary \ref{primaC} and Theorem \ref{thm:H_2} for a quick
view of their properties.

\subsection{Curve $\sH_1$ parameterizing marked lines}
\label{subsection:H1}~

\subsubsection{Construction of $\sH_1$ and marked lines}~

Set $\sH_1:=\varphi^{-1} C\subset \mP$ and $M:=M_d$.
We begin with a few corollaries of Proposition \ref{prop:Cd}:
   
\begin{cor}\label{primaC}
If $d\geq 2$, then 
 $\sH_1$ is a smooth curve of genus $d-2$
with the triple cover $\sH_1\to C$.
In particular, if $d\geq 5$, then 
$\sH_1$ is a smooth non-hyperelliptic trigonal curve of genus $d-2$.
\end{cor}

\begin{proof}  
By Propositions \ref{prop:FN} (1) and \ref{prop:Cd} (1),
it holds that $\sH_1$ is smooth and
the ramification for $\sH_1\to C$ is simple.
Since $B_{\varphi}\in |-K_{B}|$ and $d={\rm{deg\,{C}}}$,
we can compute $g(\sH_1)$ by the Hurwitz formula:
\[
  \text{$2g(\sH_1)-2=3\times (-2)+d\times 2$, equivalently,  
  $g(\sH_1)=d-2$}.      
\]
\end{proof}

\begin{cor}
    \label{orabisecanti}
The number of nodes of $M$ is $s:=\frac{(d-2)(d-3)}{2}$, whence
$C$ has $\frac{(d-2)(d-3)}{2}$ bi-secant lines on $B$.
\end{cor} 

\begin{proof}    
By Proposition \ref{prop:Cd} $(3)$,
we see that $\pi_{|\sH_1}\colon\sH_1\rightarrow M$ is birational
and $p_a(M)=\frac{(d-1)(d-2)}{2}$. Then by $g(\sH_1)=d-2$,
the number of nodes of $M$ is $\frac{(d-1)(d-2)}{2}-(d-2)=\frac{(d-2)(d-3)}{2}$. 
The latter half follows since a bi-secant line of 
$C$ corresponds to a node of $M$. 
\end{proof}

Now we select some lines on $B$ which we use in the sequel. Note that
\[
\sH_1=\{([l],t)\mid [l]\in M, t \in C\cap l\} \subset M\times C,
\]
and the elements of $\sH_1$ deserve a name:

\begin{defn}
A pair of a secant line $l$ of $C$ on $B$ and 
a point $t\in C\cap  l$ is called a {\it{marked line}}.

Let $(l,t)$ be a marked line.
If $C\cap l$ is one point, then 
$\{t\}=C\cap l$ is uniquely determined. 
For a bi-secant line $\beta_i$ of $C$, 
there are two choices of $t$.  
Thus $\sH_1$ parameterizes marked lines.
\end{defn}

%\medskip

%%%%%%%%%%%%%%%%%%%%%%%%%%%%%%%%%%%%%%%%%%%%%%%%%%%%%%%%%%%%%%%%%%%%%%%%%%%%%%%%%%%%%%%%%%%%%%%%%%%%%%%%%%%%%%%%%%%%%%%%%%%%%%%%%%%%%%%%%%%%%%%%%%%%%%%%%%%%%%%%

%%%%%%%%%%%%%%%%%%%%%%%%%%%%%%%%%%%%%%%%%%%%%%%%%%%%%%%%%%%%%%%%%%%%%%%%%%%%%%%%
\subsubsection{Lines on the blow-up $A$ of $B$ along $C_d$}
~
\label{subsubsection:lineA}

We prove that each marked line corresponds to a curve 
of anticanonical degree $1$ on the blow-up
$A$ of $B$ along $C$. This gives us a suitable notion
of line on $A$.

\begin{nota}
~
\begin{enumerate}[$(1)$] 
\item \parbox[t]{0.80\textwidth}{
Let $f\colon A\to B$ be the blowing up along $C$ and 
$E_{C}$ the $f$-exceptional divisor,} 

\item
 $\{p_{i1}, p_{i2}\}=C\cap \beta_{i}\subset B$,
\item
$\zeta_{ij}=f^{-1}(p_{ij})\subset E_{C}\subset A$, and
\item
$\beta'_{i}\cap \zeta_{ij}= p'_{ij}\in E_{C}\subset A$,
\end{enumerate}
where $i=1,\ldots, s$ and $j=1,2$.
\end{nota}

\begin{defn}
    \label{linea}
    We say that a connected curve $l\subset {A}$ 
    is a {\em{line}} on ${A}$ if
$-K_{{A}} \cdot l=1$ and
${{E_{C}}}\cdot l=1$.
\end{defn}

We point out that since $-K_A=f^*(-K_{B})-E_{C}$ and $E_{C}\cdot
l=1$ then
$f(l)$ is a line on $B$ {\it{intersecting}} $C$.
More precisely:
\begin{prop}\label{lineeA}
A line $l$ on $A$ is one of the following curves on $A:$
\renewcommand{\labelenumi}{\textup{(\roman{enumi})}}
\begin{enumerate}
\item
the strict transform of a uni-secant line of $C$ on $B$,
or
\item
the union $l_{ij}=\beta'_{i}\cup \zeta_{ij}$, where $i=1,\ldots, s$
and $j=1,2$.
\end{enumerate}
In particular $l$ is reduced and $p_a(l)=0$.
\end{prop}

\begin{nota}
For a line $l$ on $A$, we usually denote by $\overline{l}$ its image on
$B$.
\end{nota}

\begin{cor}\label{parameter} 
    The curve $\sH_1\subset\mP$ is the Hilbert scheme of
     the lines of $A$.
\end{cor}

\begin{proof} 
Let $\sH'_1$ be the Hilbert scheme of lines on $A$,
which is a locally closed subset of the Hilbert scheme
of $A$.
By the obstruction calculation of the normal bundles of the
components of
lines on $A$, it is easy to see that $\sH'_1$ is a smooth curve.
Denote by $\sU_1\to \sH'_1$ the universal family of the lines on 
$A$ and let $\overline \sU_1$ be the image of $\sU_1$
on $B\times \sH'_1$ (with induced reduced structure).

\begin{cla} 
\label{cla:basic}
$\overline {\sU}\to \sH'_{1}$ is a $\mP^1$-bundle.
\end{cla}
\begin{proof}[Proof of the claim]
Let $\sL$ be the pull-back of 
the ample generator of $\Pic B$
by 
\[
{\sU_{1}}\hookrightarrow 
A \times \sH'_{1}\to 
B\times \sH'_{1}\to B.
\]
Since $\varrho\colon \sU_{1}\to \sH'_{1}$ is flat
and $h^0(l, \sL_{|l})=2$ for a line $l$ on ${B}$,
$\sE:=\varrho_*\sL$ is a locally free sheaf of rank two.
$\mP(\sE)$ is nothing but the $\mP^1$-bundle
contained in $B\times \sH'_{1}$ whose fiber
is the image of a line on ${A}$.
This implies that $\mP(\sE)=\overline \sU$ as schemes and 
$\overline \sU$ is a $\mP^1$-bundle.
\end{proof}
By the claim, we have a natural morphism
$\sH'_1\to \mP^2$, 
whose image is $M$.
By Proposition \ref{lineeA} $\sH'_1 \to M$ is birational and surjective. 
Since $\sH'_1$ and $\sH_1$ are smooth, they are both normalizations of
$M$, hence $\sH'_1\simeq \sH_1$.
\end{proof}

\begin{rem}
For a bi-secant line $\beta_i$,
we have two choices of marking, $p_{i1}$ or $p_{i2}$.
We describe which line on $A$ corresponds to $(\beta_i,p_{ij})$. 
Denote by $\sU_1\to \sH_1$ the universal family of the lines on 
$A$ and consider the following diagram:
\begin{equation*}
\xymatrix{{\sU_1}\ar[d] \subset & 
A\times \sH_1 \ar[d]\\
{\overline{\sU}_1}  \subset &
B \times \sH_1.}
\end{equation*}
Then
$\sU_1\to \overline{\sU}_1$ is the blow-up
along $(C\times \sH_1)\cap \overline{\sU}_1$,
which is the union of a section of $\overline{\sU}_1\to \sH_1$
consisting markings
and finite set of points
$(p_{i,3-j}, [\beta_i, p_{ij}])$.
Thus 
the marked line $(\beta_i, p_{ij})$ corresponds to
the line $l_{i,3-j}$.
\end{rem}

%%%%%%%%%%%%%%%%%%%%%%%%%%%%%%%%%%%%%%%%%%%%%%%%%%%%%%%%%%%%%%%%%%%%%%%%%%%%%%%%
%%%%%%%%%%%%%%%%%%%%%%%%%%%%%%%%%%%%%%%%%%%%%%%%%%%%%%%%%%%%%%%%%%%%%%%%%%%%%%%%%%%%%%%%%%%%%%%%%%%%%%%%%%%%%%%%%%%%%%%%%%%%%%%%%%%%%%%%%%%%%%%%%%%%%%%%%%%%%%%%
\subsection{Surface $\sH_2$ parameterizing marked conics}
\label{subsection:H2}
~

Now we define a notion of {\it{conic}} on $A$. We proceed as in the
case of lines, first defining the notion of {\it{marked conic}}.

%%%%%%%%%%%%%%COSTRUZIONE di H_{2}%%%%%%%%%%%%%%%%%%%%%%%
\subsubsection{Construction of $\sH_2$ and marked conics}
\label{subsubsection:H_2}

\begin{defn}\label{conichemarked}
A pair of a multi-secant conic $q$ on $B$ 
and a zero-dimensional subscheme $\eta \subset C$
of length two contained in $q_{|C}$ 
is called a {\it{marked conic}}.
%We call the union of two marked lines which is a marked conic
%{\it{a marked line pair}}.\footnote{Needed ?}
\end{defn}

From now on, we assume that $d\geq 3$.

Marked conics are parameterized by
\[
\sH'_2:=\{([q],[\eta])\mid [q]\in \overline{\sH}'_2, \eta\subset 
q_{|C}\}
\subset \overline{\sH}'_2\times S^2 C
\]
with reduced structure,
where
$\overline{\sH}'_2\subset \mP^4$ is the locus of
multi-secant conics of $C$ on $B$.

By Corollary \ref{cor:uniqueconic} and $d\not =1$,
the natural projection of $\sH'_2\to S^2 C$ 
is one to one outside $[\beta_{i|C}]$ and the diagonal of $S^2 C$.
%thus 
%by the Zariski main theorem, it is an isomorphism
%outside $[\beta_{i|C}]$ and the diagonal of $S^2 C$.

We denote by 
$e'_i$ the fiber of  
$\sH'_2\to S^2 C$ over a $[\beta_{i|C}]$.
Since $B$ is the intersection of quadrics,
any conic cannot intersect a line twice properly.
Thus any conic $\supset \beta_{i|C}$ contains $\beta_i$.
This implies that $e'_i\simeq \mP^1$, and $e'_i$
parameterizes marked conics of the form
\[
\{([\beta_i\cup \alpha], [\beta_{i|C}])\mid 
\text{$\alpha$ is a line such that $\alpha\cap \beta_i\not
=\emptyset$}\}.
\]

Over the diagonal of $S^2 C$,
$\sH'_2\to S^2 C$ is finite since 
for $t\in C$,
there exist a finite number of reducible conics with $t$
as a singular point or conics tangent to $C$ at $t$.

Hence $\sH'_2$ is the union of 
the unique two-dimensional component, which dominates $S^2 C$,
and possibly lower dimensional components mapped into
the diagonal of $S^2 C$ or $e'_i$.
%Moreover $\sH'_2\to S^2 C$ is 
%isomorphic outside $[\beta_{i|C}]$. 
Note that $\sH'_2\to \overline{\sH}'_2$ is finite since
choices of markings of a multi-secant conic of $C$ is finitely many
by $d\geq 3$. 
\begin{cla}
$e'_i$ is contained 
in the unique two-dimensional component of $\sH'_2$.
\end{cla}

\begin{proof}
We have only to prove that $\overline{\sH}'_2$ is two-dimensional
near the generic point of the image of $e'_i$ since $\sH'_2\to \overline{\sH}'_2$ is one to one near the generic point of the image of $e'_i$.
Let $\sV_2\to \sH_2^{B}\simeq \mP^4$ be the universal family of
conics on $B$ and $\overline{\sH}''_2$
the inverse image of $C\times C$ 
by $\sV_2\times _{\mP^4} \sV_2 \to B\times B$.
%and $\overline{\sH}''_2\to \overline{\sH}'_2$ is finite of degree two.
Since the morphism $\sV_2\times _{\mP^4} \sV_2 \to \sV_2\to \mP^4$ is flat,
$\sV_2\times _{\mP^4} \sV_2$ is purely six-dimensional.
Thus any component of $\overline{\sH}''_2$ 
has dimension greater than or equal to two.
Though the inverse image of
the diagonal of $C\times C$ is three-dimensional,
any other component of $\overline{\sH}''_2$ is at most two-dimensional 
by a similar investigation to $\sH'_2$.
Thus $\overline{\sH}'_2$ is two-dimensional
near the generic point of the image of $e'_i$ since
$\overline{\sH}'_2$ is the image of the two-dimensional part of 
$\overline{\sH}''_2$ by $\sV_2\times _{\mP^4} \sV_2\to \mP^4$
near the generic point of the image of $e'_i$.
\end{proof}

\begin{nota}
Let $\sH_2$ be the normalization of the unique two-dimensional component of
$\sH'_2$ and $\overline{\sH}_2\subset \overline{\sH}'_2$ 
the image of $\sH_2$. 
Denote by $\eta$ the natural morphism $\sH_2\to S^2 C$.
Set 
\[
c_i:=[\beta_{i|C}] \in S^2 C\simeq \mP^2,
\]
and
\[ 
e_i:={\eta}^{-1}(c_i),
\] 
where $i=1,\ldots, s$.
\end{nota}

By the above consideration,
$\eta\colon \sH_2\to S^2 C$ is isomorphic outside $[\beta_{i|C}]$
by the Zariski main theorem,
and $\sH_2\to \overline{\sH}_2$ is the normalization.
Thus we see that $\sH_2$ parameterizes marked conics
in one to one way
outside the inverse image of $c_i$. We need to understand the inverse
image by $\eta$ of the diagonal.

\begin{cla}
\label{cla:reduced}
Assume that $([q],[2b])\in \sH_2$ for $b\in C$ and a conic $q$. Then
\begin{enumerate}[$(1)$]
\item
$q$ is reduced,
\item
if $q$ is smooth at $b$, then $q$ is tangent to $C$ at $b$, and
\item
if $q$ is singular at $b$, then
the strict transform of $q$ is connected on $A$.
Moreover, $b\not \in \beta_i$ nor $B_{\varphi}$.
\end{enumerate}
\end{cla}

\begin{proof}
We use the double projection from $b$.
By Proposition \ref{prop:proj2} (4) and a degeneration argument, 
$q$ corresponds to the fiber of $\pi_{2b}$ through 
the point $t'$ in $C''_b\cap E'_b$ coming from $t:=C'_b\cap E_b$.

(1) Assume by contradiction that $q$ is non-reduced.
By Proposition \ref{prop:conics}, $q$ is a multiple of 
a special line $l$. By Proposition \ref{prop:Cd1} (4),
$l$ is a uni-secant line of $C$.
Let $m$ be the other line through $b$ (by generality of $C$, we have $l\not =m$). Let $l'$ and $m'$ be the strict transforms of $l$ and $m$ on $B_b$ respectively. By Proposition \ref{prop:proj2} (4),
the fiber of $\pi_{2b}$ through $t'$ is the strict transform
of the line in $E_b$ joining $l'\cap E_b$ and $m'\cap E_b$.
Hence by the assumption,
$l'\cap E_b$, $m'\cap E_b$ and $C'_b\cap E_b$ are collinear.
By dimension count similar to 
the proof of Proposition \ref{prop:Cd1}, 
we can prove that a general $C$ does not satisfy this condition.

(2) This follows from the previous discussion.

(3) Set $q=l_1\cup l_2$, where
$l_1$ and $l_2$ are the irreducible components of $q$,
and let $l'_i$ be the strict transform of $l_i$ on $B_b$.
By (1), it holds $l_1\neq l_2$.
Then the fiber of $\pi_{2b}$ corresponding to $q$ 
is the strict transform of the line on $E_b$ through $E_b\cap l'_1$
and $E_b\cap l'_2$.
Note that $A$ is obtained from $B_b$ by blow- up $B_b$ along $C'_b$
and then contracting the strict transform of $E_b$.
Thus the former half of the assertion follows.
The latter half follows again by simple dimension count.
\end{proof}

%%%%%%%%%%%%%%COSTRUZIONE di H_{2}:FINE %%%%%%%%%%%%%%%%%%%%%%%

%&&&&&&&&&&&&&&&CONICS ON A&&&&&&&&&&&&&&&&&&&&&&&&&&&&&&&&&&&&&&&&
\subsubsection{Conics on $A$}
\label{subsubsection:coniche}

\begin{defn}
    \label{conica}
    We say that a connected and reduced curve $q\subset {A}$ 
    is a {\em{conic}} on ${A}$ if
    $-K_{{A}} \cdot q = 2$ and ${{E_{C}}}\cdot q=2$.
\end{defn}

Using this definition, we can classify conics on $A$
similarly to Proposition \ref{lineeA}:

\begin{prop}\label{classificazioneconicheinA}
Let $q$ be a conic on ${A}$.
Then $\overline q:=f(q)\subset B$ is a multi-secant conic of 
$C$.
Moreover one of the following holds$:$
\renewcommand{\labelenumi}{\textup{(\alph{enumi})}}
\begin{enumerate}
\item
$\overline q$ is smooth at $\overline q\cap C$. 
$q$ is the union of the strict transform $q'$ 
of $\overline q$
and $k-2$ distinct fibers $\zeta_1,\dots, \zeta_{k-2}$ of $E_C$ 
such that $\zeta_i\cap q' \not =\emptyset$,
\item
$\overline{q}$ is the union of two uni-secant lines $\overline l$ and
$\overline m$ such that 
$C\cap \overline l\cap \overline m \not =\emptyset$.
$q$ is the union of the strict transforms $l$ and $m$ 
of $\overline l$ and
$\overline m$ respectively $($we assume that $l\cap m\not =\emptyset$$)$, or 
\item
$\overline{q}$ is the union of $\beta_i$ and a line $\overline{r}$ through a
$p_{ij}$. $q$ is the union of the fiber $\zeta_{ij}$ over $p_{ij}$
and the strict transforms $\beta'_i$ and $r'$ 
of $\beta_i$ and $\overline{r}$ respectively.
%\footnote{Even if the strict transforms $\overline q'$ and $r'$ 
%of $\overline q$ and $r'$, no problem ? We had better disprove this
%for the proof of Prop. 6.12 ?}
\end{enumerate}  
\end{prop}

%\begin{proof} 
%By Definition \ref{conica} (i), (ii), and (iii)
%and the formula $-K_A=f^* (-K_{B})-E_C$,
%we see that $\overline q$ is a conic.
%Thus 
%$q$ is the union of the strict transform $q'$ of $\overline{q}$
%and some fibers of $E_C\to C$.
%By the conditions of Definition \ref{conica}, 
%we have the conclusion.
%\end{proof}

%%%%%%%%%%CONICHE=CONICHE MARCATE%%%%%%%%%%%%%%

\begin{nota}
We usually denote by $\overline{q}\subset B$ the image of a conic $q$ on $A$.
\end{nota}
 
Let $\sH^A_2$ be the normalization of the two-dimensional part of the Hilbert scheme of conics on $A$,
which is a locally closed subset of the Hilbert scheme
of $A$.
Let
$\mu \colon {\sU_2}\to \sH^A_2$ be the pull-back of
the universal family of conics on ${A}$.
% clearly some fibers
%in ${\sU_2}\subset A\times \sH^A_2 $ are not conic 
%in the standard sense. 

\begin{lem}\label{useful?}
    Let
$\overline{\sU}_2$ be the image of ${\sU_2}$
on $B\times \sH^A_2$ $($with induced reduced structure$)$ then
$\overline \sU_2\to \sH^A_2$ is a conic bundle. 
%In particular $\overline{\sU}_2$ is normal. 
\end{lem}

\begin{proof}
The proof is similar to that of Claim \ref{cla:basic}.

Let $\sL$ be the pull-back of 
the ample generator of $\Pic B$
by 
\[
{\sU_{2}}\hookrightarrow 
A \times \sH^A_2\to 
B\times \sH^A_2\to B.
\]

Since $\mu\colon {\sU_2}\to \sH^A_2$ is flat
and $h^0(q, \sL_{|q})=3$ for a conic
$q$ on ${A}$ (recall that $q$ is reduced), then
$\sE:=\mu_*\sL$ is a locally free sheaf of rank $3$.
Letting
$\mP^{6}=\langle B \rangle$, $\mP(\sE)$ is the $\mP^2$-bundle
contained in $\mP^{6}\times\sH^A_{2}$ whose fiber
is the plane spanned by the image of a conic on ${A}$.
Let $\mathcal{Q}:=(B\times \sH^A_2) \cap \mP(\sE)$,
where the intersection is taken in $\mP^{6}\times \sH^A_2$.  
A scheme theoretic fiber of $\sQ \to \sH^A_2$ 
is the image of a conic of ${A}$ since
$B$ is the intersection of quadrics.
Then $\sQ=\overline{\sU}_2$ as schemes and 
$\overline{\sU}_2$ is a conic bundle.
\end{proof}

\begin{prop}\label{equality}
There exists a natural bijection between
the set of marked conics belonging to $\sH_2$
and the set of conics on $A$.
Moreover,
the two surfaces $\sH_2^{A}$ and $\sH_2$ are isomorphic.
\end{prop}

\begin{proof} 
The first assertion follows from
Claim \ref{cla:reduced} (1) and (3), 
and Proposition \ref{classificazioneconicheinA}.

By Lemma \ref{useful?},
there exists a natural morphism $\overline{\nu}\colon \sH^A_2\to \overline{\sH}'_2$.
By Proposition \ref{classificazioneconicheinA},
$\overline{\nu}$ is finite and birational, 
hence $\overline{\nu}$ lifts to the morphism $\nu\colon \sH^A_2 \to \sH_2$
since $\sH_2\to \overline{\sH}_2$ is the normalization.
By the Zariski main theorem, $\nu$ is an inclusion.
By Claim \ref{cla:reduced} (1) and (3), 
and Proposition \ref{classificazioneconicheinA}, 
$\nu$ is also surjective.
\end{proof}
 
 By Proposition \ref{equality} we can pass freely from conics on $A$,
 that is elements of $\sH_2^{A}$ to marked conics and vice-versa
 according to the kind of argument we will need. In particular we can
 speak of the universal family $\mu\colon \sU_{2}\rightarrow \sH_{2}$ of marked 
 conics meaning $\sU_{2}:=\sU^{A}_{2}$ and $\sH_2^{A}$ identified to
 $\sH_{2}$ via $\nu$. 
 
\begin{cor}
\label{cor:H2}
The Hilbert scheme of conics on $A$ is an irreducible surface
$($and $\sH_2$ is the normalization$)$.
The normalization is injective, 
namely, $\sH_2$ parameterizes conics on $A$ in one to one way.  
\end{cor}

\begin{proof}
By Proposition \ref{classificazioneconicheinA},
the image of $\sH_2$ in the Hilbert scheme of $A$
parameterizes all the conics on $A$, thus the first part follows.

For the second part,
we have already seen that $\sH_2$
parameterizes marked conics belonging to $\sH_2$
in one to one way outside $\cup_i e_i$.
Thus, by Proposition \ref{equality},
$\sH_2$ parameterizes conics on $A$ in one to one way outside $\cup_i e_i$.
Let $\alpha$ be a general line intersecting $\beta_i$,
and $\alpha'$ the strict transform of $\alpha$ on $A$.
By easy obstruction calculation, we see that
the Hilbert scheme of conics on $A$ is smooth at $[\beta'_i\cup \alpha']$.
Thus general points of $e_i$ also parameterizes conics on $A$
in one to one way.
Then, however, since $e'_i\simeq \mP^1$, where
$e'_i$ is the inverse image of $[\beta_{i|C}]$ by $\sH'_2\to S^2 C$,
it holds that $e_i\simeq e'_i \simeq \mP^1$
($\sH_2\to S^2 C$ has only connected fibers).
This implies the assertion.
\end{proof}
\medskip

%%%%%%%%%%CONICHE=CONICHE MARCATE:FINE%%%%%%%%%%%%%%

In subsection \ref{subsubsection:H'_2},
we complete the description of $\sH_2$.
In \ref{subsubsection:qf} and \ref{lineconic},
we give some preliminary results for this purpose.

\subsubsection{Quasi-finiteness of $\psi\colon \sU_2\to A$}
\label{subsubsection:qf}

\begin{nota}
For a point $b\in C$, set
\[
L_b:=\overline{\{[q]\in \sH_2 \mid \exists b'\not =b,
f(q)\cap C=\{b, b' \} \}}.
\]
By Corollary \ref{cor:uniqueconic}, 
$\eta(L_b)$ is a line in $S^2 C \simeq \mP^2$.
\end{nota}

Let $\psi\colon\sU_{2}\to A$ be the morphism obtained
    via the universal family $\mu\colon\sU_{2}\to\sH_{2}$. 
The following result refines Proposition \ref{prop:Cd2}.
%From now on, we assume $d\geq 5$ throughout the paper
%since we need Proposition \ref{prop:Cd2}.
Here we need this result technically for the discussion
in \ref{lineconic} but
this is important for the proof of the main result and
is refined again in \ref{subsection:finite}
(Proposition \ref{finitezzafinale}).

From now on in this paper, we assume that $d\geq 5$.
\begin{prop}\label{fuoripiatto}
 The morphism $\psi$ is finite of degree $n=\frac{(d-1)(d-2)}{2}$ 
 and flat outside $\cup_{i=1}^{s} \beta'_i$.     
 \end{prop}

\begin{proof}  
    Let
    $a\in A\setminus\cup_{i=1}^{s}\beta'_{i}$ and set $b:=f(a)$.
    If $b\not \in C$, then the finiteness of $\psi$ over $a$
    follows from Corollary \ref{cor:finite}. 
    Moreover, by Corollary \ref{cor:n},
    the number of conics through a general $a$ is $n$.
    Thus $\deg \psi=n$. We will prove that $\psi$ is finite over 
    $a\in E_C\setminus \cup_{i=1}^{s} \beta'_i$.
    Once we prove this,
    the assertion follows. Indeed,
    $\sU_2$ is Cohen-Macaulay since 
    $\sH_{2}$ is smooth
    and any fiber of $\sU_2\to\sH_2$ is reduced, thus  
    $\psi$ is flat.
    
    Let $a\in E_C\setminus \cup_{i=1}^{s} \beta'_i$.
    %By Corollary \ref{cor:finite}, 
    The assertion is equivalent to that
    only finitely many conics belonging to $L_b$ pass through $a$.
    If $b\not \in \cup_{i=1}^{s} \beta_i$, then $L_b$ is irreducible.
    If $b\in \cup_{i=1}^{s} \beta_i$, then $L_b=L'_b\cup e_i$,
    where $L'_b$ is the strict transform of $\eta(L_b)$ whence is irreducible.
    Note that almost all the conics belonging to $e_i$
    does not pass through $a \not \in \cup_{i=1}^{s} \beta'_i$.
    Let $S_{b}\subset {A}$ be the locus swept by the 
    conics of the family $L_{b}$ if $b\not \in \cup_{i=1}^{s} \beta_i$, 
    or the locus swept by the 
    conics of the family $L'_{b}$ if $b\in \cup_{i=1}^{s} \beta_i$.
    Then $S_b$ is irreducible.
    Let $\overline{S}_b:=f(S_b)$,
    $\overline{S}'_b$ and $\overline{S}''_b$ 
    the strict transforms of $\overline{S}_b$ on $B_b$ and 
    $B'_b$ respectively.
    Then $\overline{S}''_b=\pi_{2b}^* C_b$. 
    %\footnote{We need to assume that $d\geq 5$, 
    %where we use Corollary \ref{cor:finite}.}
    Let $d_b:=\deg C_b$. By Proposition \ref{prop:proj2} (2),
    $d_b=d-2$ if $b\not \in \cup_{i=1}^{s} \beta_i$, or
    $d-3$ if $b\in \cup_{i=1}^{s} \beta_i$.
    Since $\overline{S}''_b\sim d_b L$ and $L=H-2E'_b$, we have
    $\overline{S}'_{b|{E_b}}$ is a curve of degree $2d_b$ in $E_b\simeq \mP^2$.
    
    Since $A$ is obtained from $B_b$ by blowing up $C'_b$ and then
    contracting the strict transform of $E_b$,
    a point $a$ over $b$ corresponds to a line $l_a$ in $E_b$ through 
    $t:=E_b\cap C'_b$. The image on $B_{b}$ of the strict transform 
    of a conic on $A$ through $a$
    intersects $E_b$ at a point of $l_a\cap \overline{S}'_b$.
    If $C''_b$ does not intersect fibers of $\pi_{2b}$ contained in  $E'_b$, 
    then
    $\overline{S}''_{b|{E'_b}}$ is irreducible. 
    Thus no $l_a$ is 
    contained in $\overline{S}'_{b|{E_b}}$ and we are done.
    Assume that  
    $C''_b$ intersects a fiber $l'$ of $\pi_{2b}$ contained in $E'_b$.
    This is a situation as in Claim \ref{cla:reduced} (3), hence
    $b\not \in B_{\varphi}$ nor $b \not \in \cup_{i=1}^{s} \beta_i$
    for a general $C$.
    %\footnote{Better prove before Claim \ref{cla:connected}.}
    Since $L_b$ is irreducible by $b \not \in \cup_{i=1}^{s} \beta_i$,
    it suffices to prove the finiteness and nonemptyness of the set of conics
    through a general point $a$ over $b$. 
    Equivalently, we have only to show that 
    a general $l_a$ intersects 
    $\overline{S}'_{b|E_b}$ outside $t$.    
    Since $l'$ intersects $C''_b$ simply at one point,
    $C_b$ is smooth at the image $t'$ of $l'$ on $\mP^2$.     
    Thus   
    $\overline{S}'_{b|{E_b}}=C'''_b+l$, 
    where $C'''_b$ and $l$ are the strict transforms
    of $C_b$ and $l'$. Note that $C'''_b$ is smooth at $t$ and
    $\deg C'''_b=2d_b-1=2d-5\geq 5$ by $d\geq 5$.
    %\footnote{$d\geq 4$ ?}
    Thus a general $l_a$ intersect $C'''_b$ outside $t$.
\end{proof}

\begin{rem}
Though we do not need it later,
we describe the fiber of $\psi$ 
over a general point $a \in E_C\setminus \cup_{i=1}^{s} \beta'_i$
for reader's convenience.

Set $b:=f(a)$. 
As in the proof of Proposition \ref{fuoripiatto}, 
a point $a$ over $b$ corresponds to a line $l_a$ in $E_b$
passing through $E_b\cap C'_b$.
By Lemma \ref{lem:Cb},
it holds that $\deg C_b=d-2$ and $C_b$ has $\frac{(d-3)(d-4)}{2}$ simple nodes
for a general $b\in C$. This means that  
$\frac{(d-3)(d-4)}{2}$ tri-secant conics pass through $b$. 
By Proposition \ref{classificazioneconicheinA}, corresponding to a tri-secant conic
$\overline q$,
there is a unique conic $q$ on $A$ containing the fiber of $E_C$
over $b$ and such a conic on $A$ contains $a$. 
Thus we obtain $\frac{(d-3)(d-4)}{2}$ conics through $a$.
By definition of $L_b$, these conics do not belong to $L_b$.

We need more
$n-\frac{(d-3)(d-4)}{2}=2d-5$ conics through $a$.
We show that there exist $2(d-2)-1$ conics through $a$ on $A$
coming from the family parameterized by $L_b$.
We use the notation of the proof of Proposition \ref{fuoripiatto}.
For a general $b\in C$,
$C''_b$ does not intersect fibers of $\pi_{2b}$ contained in $E'_b$. 
Thus $\overline{S}'_{b|E_b}$ is an irreducible curve
of degree $2(d-2)$ on $E_b$.
% by $L=H-2E'_b$ and $\overline{S}''_b\sim d_b L$. 
Thus there are $2(d-2)$ intersection points of
$\overline{S}'_{b|E_b}$ and $l_a$. Among those, the intersection
point $C'_{b}\cap E_b$ does not correspond to a conic on $A$
through $a$ since it comes from the tangent of $C$. 
Thus we have $2(d-2)-1$ conics as desired.
\end{rem}

\subsubsection{Intersection of lines and conics on $A$}
\label{lineconic}~

To understand better $\eta\colon \sH_2\to\mP^2$ 
we need to find special loci inside $\sH_{2}$. A natural step
is to study the locus of conics which intersect a fixed line. 
This locus turn out to be a good divisor of $\sH_2$.

Let ${\sU}'_1 \subset {\sU_2}\times \sH_1$ be the
pull-back of ${\sU_1}$ via the following diagram:
\begin{equation}
\label{eq:familyB1}
\xymatrix{\sU'_1\subset {\sU_2}\times \sH_1 \ar[d] 
\ar[r] & A \times \sH_1 \supset \sU_1 \ar[d]\\
\widehat{\sD}_1 \subset \sH_2\times \sH_1 \ar[r] & \sH_1,} 
\end{equation}
where $\widehat{\sD}_1$ is the image of $\sU'_1$ on 
$\sH_2\times \sH_1$.

By definition
\[
\widehat{\sD}_1=\{([q], [l])\mid q \cap l \not =\emptyset\}
\subset \sH_2\times \sH_1.
\]
%\noindent where the intersection is between a conic and a line on $A$.

First we need to know which component of $\widehat{\sD}_1$
is divisorial or dominates $\sH_1$.
For this purpose,
we study mutual intersection of a conic and a line in special cases. 
Let $\sF\subset\sH_{2}\times\sH_{1}$ be
    the image in $\sH_2\times \sH_1$ of 
the inverse image of $((\cup \beta'_i)\times \sH_1) \cap \sU_1$; that 
is
\[
\sF:=\{([q],[l])\mid q \cap \beta'_i\cap l\not =\emptyset\}. 
\]
A point $([q],[l])\in\sF$ iff
i) $l=l_{ij}(:=\beta'_i\cup \zeta_{ij})$ and 
$q\cap \beta'_i\not =\emptyset$ or ii)
$l\not =l_{ij}$, and 
$q \cap \beta'_i\cap l\not =\emptyset$. 
For every $i=1,\ldots ,s$, $j=1,2$ the family of those $([q],[l])$
which satisfies i) or ii) has dimension one and clearly does not dominate
$\sH_{1}$. 

\begin{cor}\label{divisoriali}
Any component of $\widehat{\sD}_1$ which is not contained in
$\sF$ dominates $\sH_1$.
Moreover, any non-divisorial component of $\widehat{\sD}_1$
outside $\sF$
$(\text{if it exists})$ is
a one-dimensional component whose generic point
parameterizes reducible conics, namely,
a one-dimensional component of 
\[
\{([q],[l])\mid l\subset q\}.
\]
%Except possibly the components of $\sF$
%and some subschemes of dimension at most one which
%do not dominate $\sH_1$,
%any other component of $\widehat{\sD}_1$ is 
%divisorial and dominates $\sH_{1}$.
\end{cor}

\begin{rem}
Here we leave the possibility that
a one-dimensional component whose generic point
parameterizes reducible conics is contained in 
a divisorial component of $\widehat{\sD}_1$.
We, however, prove that 
this is not the case in Corollary \ref{noncontenere}.
Hence, finally, the fiber of $\widehat{\sD}_1\to \sH_1$ over
a general $[l]\in \sH_{1}$ parameterizes 
conics which properly intersect $l$. 
\end{rem}

\begin{proof} 
    By Proposition \ref{fuoripiatto}, $\sU_2\to A$ is finite and flat outside $\cup
    \beta'_i$. Then $\sU_{2}\times\sH_{1}\to A\times\sH_{1}$ is flat
    outside $(\cup \beta'_i)\times \sH_1$. By base change,
    $\sU'_{1}\rightarrow \sU_{1}$ is flat and finite outside 
    $((\cup \beta'_i)\times \sH_1) \cap \sU_1$. Then every 
    irreducible component of $\sU'_1$ which is not mapped to 
$((\cup \beta'_i)\times \sH_1) \cap \sU_1$ is two-dimensional, and
dominates $\sU_1$, hence dominates $\sH_1$.
Therefore
any component of $\widehat{\sD}_1$ which is not contained in
$\sF$ dominates $\sH_1$.
%The subscheme $[(\cup \beta'_i)\times \sH_1] \cap \sU_1$
%does not dominate $\sH_1$.

We find a possible non-divisorial component of
$\widehat{\sD}_1$ outside $\sF$.
Let $\gamma\subset \sU'_1$ be a curve mapped to a point, say, 
$([q],[l])$ 
on $\sH_2\times \sH_1$. The image of $\gamma$
on $A$ is an irreducible component of $q$, say, $q_1$. 
The image of $\gamma$ on $\sU_1$ 
is $q_1 \times [l]$, thus $q_1$ is also an irreducible
component of $l$.
We have the following three possibilities:

\begin{enumerate}[$(1)$]
\item
$l$ is irreducible, hence $q_1=l$ and
$q=l\cup m$, where $m$ is another line. 
Such $([q],[l])$ form the one-dimensional family of reducible conics,
%\footnote{remove marked line pair.}
\item
$l=l_{ij}$ and $\beta'_i\subset q$.
Namely $[q]\in e_i$, or $q=\beta'_i\cup \alpha\cup \zeta_{ik}$,
where $\alpha$ is the strict transform of a line on $B$
intersecting $\beta_i$ and $C$ outside $\beta_i\cap C$, or
\item
$l=l_{ij}$ and $\zeta_{ij}\subset q$ and
$f(q)$ is a tri- or quadri-secant conic of $C$
such that $p_{ij}\in f(q)$.
\end{enumerate}
Thus we have the second assertion.
\end{proof}

\begin{nota}
Let ${\sD}_1 \subset \sH_2\times \sH_1$ be 
the divisorial part of $\widehat{\sD}_1$. 
Since $\sH_1$ is a smooth curve ${\sD}_1 \to \sH_1$ is flat.  
Let ${D}_l$ be the fiber of 
${\sD}_1\to \sH_1$ over $[l]\in \sH_1$. Clearly we can write
$D_{l}\hookrightarrow \sH_{2}$.
\end{nota}

%Next two lemmas are basic to understand the geometry of $\sH_{2}$.

\subsubsection{Description of $\sH_2$}
\label{subsubsection:H'_2}~

Now we reach the precise description of $\sH_2$.

%In the proof of the following theorem,
%we investigate the intersection
%of several $D_l$, and $L_b$.
%By using Lemmas \ref{lem:twolines},
%\ref{lem:twolines'} and Corollary \ref{cor:n},
%we can easily compute the number of points contained in
%the intersections but
%we need to show the intersection is simple
%to compute the intersection numbers.
%This is the 

\begin{thm}
\label{thm:H_2}
\begin{enumerate}[$(1)$]
\item
The morphism $\eta\colon \sH_2\to \mP^2$
is the blow-up at $c_1,\dots, c_s$
and $e_i$ are $\eta$-exceptional curves. 
It holds:
\[
D_l\sim (d-3)h-\sum_{i=1}^{s} e_i,
\]
where $h$ is the strict transform of a general line on $\mP^2$.
\item
\[
h^1(\sH_{2},\sO_{\sH_{2}}((d-4)h-\sum_{i=1}^{s} e_i))=0.
\]
\item
$|D_l|$ is base point free.
In case of $d=5$, the image of $\Phi_{|D_l|}$ is $\check{\mP}^2$.
In case of $d\geq 6$, $D_l$ is very ample and
$|D_l|$ embeds $\sH_2$ into $\check{\mP}^{d-3}$.

Here we use the dual notation $\check{\mP}^{d-3}$
for later convenience.

\item
If $d\geq 6$, then 
$\sH_2\subset \check{\mP}^{d-3}$ is 
projectively Cohen-Macaulay,
equivalently,
\[
\text{$h^i(\check{\mP}^{d-3}, \sI_{\sH_2}(j))=0$
for $i=1,2$ and $j\in \mZ$,}
\]
where $\sI_{\sH_2}$ is the ideal sheaf of $\sH_2$ in
$\check{\mP}^{d-3}$. 
Moreover, $\sH_2$ is the intersection of cubics.
 
\end{enumerate}
\end{thm}

\begin{rem}
If $d\geq 6$, then $\sH_2\subset \check{\mP}^{d-3}$ is so called the 
{\em{White surface}} (see \cite{White} and \cite{Gimi}).
In \cite{Man},
the White surface attains the maximal degree among 
projectively Cohen-Macaulay rational surfaces in
a fixed projective space.
\end{rem}

\begin{proof}
(1) 
First we compute the intersection number $D_l\cdot L_b$
for general $l$ and $b$
(this intersection number will be well-defined
since the intersection points of $D_l$ and $L_b$
are contained in the smooth locus of $\sH_2$).  
We prove that $D_l$ and $L_b$ intersect simply.
Indeed,
let $\pi_C\colon C\times C\to S^2 C$ be the natural projection
and $L'_b$ a ruling of $C\times C\to C$ in one fixed direction 
such that $\pi_C(L'_b)=\eta(L_b)$.
By applying the Bertini theorem to $|L'_b|$, 
we see that $\pi_C^*\eta(D_l)$ and $L'_b$ intersect simply for a general
$b\in C$ whence $\eta(D_l)$ intersects $\eta(L_{b})$ simply since 
$\pi_C$ is \'etale at $\pi_C^*\eta(D_l)\cap L'_b$. Then $D_l$ intersects $L_b$ 
simply since $\eta$ is isomorphic at $D_{l}\cap L_{b}$.
Thus we have only to count the number of points in $D_l\cap L_b$,
which is $d-3$ by Corollary \ref{cor:n} (3).
Now we see $D_l\cdot L_b=d-3$ whence 
$\eta(D_l)$ is a curve of degree $d-3$.

Second, we compute the intersection number $D_{l_1}\cdot D_{l_2}$
for two general lines $l_1$ and $l_2$ on $A$.
The images $\overline{l}_1:=f(l_1)$ and $\overline{l}_2:=f(l_2)$ be 
two general secant lines of $C$
such that ${\overline{l}_1}\cap {\overline{l}_2}=\emptyset$.
By Lemma \ref{lem:twolines},
$\# (D_{l_1}\cap D_{l_2})=\frac{(d-3)(d-4)}{2}$.
This immediately gives for the intersection product
$D_{l_1}\cdot D_{l_2}\geq \frac{(d-3)(d-4)}{2}$.
Unfortunately, we cannot show the intersection is simple apriori
so we need some argument. 
On the other hand, $D_l\cap e_i\not =\emptyset$ for a general $l$
since $D_l\cap e_i$ contains the point corresponding to
a marked conic
$(\beta_i\cup \alpha, \beta_{i|C})$, where $\alpha$ is the unique line
intersecting $\beta_i$ and $l$.
Moreover, for two general $l_1$ and $l_2$, 
$D_{l_1}\cap D_{l_2}\cap e_i=\emptyset$, and
$D_{l_1}\cap e_i$ and $D_{l_2}\cap e_i$ are contained in the smooth locus
of $\sH_2$. Thus, by taking the minimal resolution
of $\sH_2$ near $e_i$ if necessarily, we can see that
$D_{l_1}\cdot D_{l_2}\leq
(d-3)^2-s=\frac{(d-3)(d-4)}{2}$.
Therefore
$D_{l_1}\cdot D_{l_2}=\frac{(d-3)(d-4)}{2}$. 
Moreover $e_{i}^{2}=-1$ and since 
$e_{i}\cap e_{j}=\emptyset$ we obtain that $\eta\colon \sH_2\to \mP^2$
is the blow-up at $c_1,\dots, c_s$.
Consequently,
$D_l\sim (d-3)h-\sum_{i=1}^{s} e_i$ for a general $[l] \in \sH_1$, and,
by the flatness of $\sD_1 \to \sH_1$,
that holds for any $[l] \in \sH_1$. 
%Let $\sD'\subset S^2 C \times \sH_1$ be the family of $\eta(D_l)$
%($[l]\in %\sH_1$).
%Since $\eta(D_l)$ is a divisor of degree $d-3$, $\sD'\to \sH_1$ is
%flat.
%Since it is easy to see that $D_l$ does not contain $e_1,\dots,e_s$,
%the strict transform $\sD$ of $\sD'$ on $\sH_2 \times \sH_1$ 
%is nothing but the family of $D_l$. Since $\sD\to \sH_1$ is flat,
%$D_l$ is linearly equivalent to $(d-3)h-\sum_{i=1}^{s} e_i$ for any
%$l$.
%Now we prove that $h^1(\sH_{2}, \sO_{\sH_{2}}(d-4)-\sum_{i=1}^{s}
%e_i))=0$. 

(2) Let $L'_{p_{ij}}=L_{p_{ij}}-e_i$ (note that $e_i\subset L_{p_{ij}}$).
We see that
$L'_{p_{ij}}\subset D_{l_{ij}}$ and 
$D_{l_{i1}}-L'_{p_{i1}}=D_{l_{i2}}-L'_{p_{i2}}$,
which we denote by $D_{\beta_i}$. 
%We have $L_{p_{ij}}\sim h-\delta_i$ since $\o(L_{s_{ij}})$ is a line
%and
%$L_{s_{ij}}\cap \delta_i$ consists of one point.
Note that
\[
D_{\beta_i}\sim (d-4)h-\sum_{k\not =i} e_k.
\]
It is easy to see that $D_{\beta_i}$ have the following properties:
\begin{eqnarray}
D_{\beta_i}\cap e_i=\emptyset.
\label{eqnarray:1}\\
D_{\beta_i}\cap D_{\beta_j}\cap D_{\beta_k}=\emptyset.
\label{eqnarray:2}
\end{eqnarray}
We only prove (\ref{eqnarray:1}).
Since $D_{\beta_i}\cap e_i\not =\emptyset$ would imply that $e_i$ is a component of
$D_{\beta_i}$, it suffices to prove that, for a general $l$,
$D_{\beta_i}\cap D_l$ does not contain
a point of $e_i$.
By Lemma \ref{lem:twolines'},
$D_{\beta_i}\cap D_l$ contains
$\frac{(d-4)(d-5)}{2}$ points
corresponding to bi-secant conics intersecting
$\beta_i$ and $l$ except conics containing $\beta_i$.
On the other hand,
we have $D_l\cdot D_{\beta_i}=\frac{(d-4)(d-5)}{2}$,
thus the conics we count in Lemma \ref{lem:twolines'}
correspond to all the intersection of $D_{\beta_i}\cap D_l$.
Consequently, $D_{\beta_i}\cap D_l$ does not contain
a point of $e_i$.

By (\ref{eqnarray:1}) and the trivial equality
\[
(d-4)h-\sum_{i\geq k+1} e_i=D_{\beta_k}+ e_1+\cdots+ e_{k-1},
\]
\noindent 
we obtain $e_k\not \subset \Bs |(d-4)h-\sum_{i\geq k+1} e_i|$.

Since $\sO_{\sH_{2}}((d-4)h-\sum_{i\geq k+1}
e_i)\otimes_{\sO_{\sH_{2}}}\sO_{e_{k}}\simeq \sO_{e_{k}}$ we have that
\[
H^0(\sH_{2},\sO_{\sH_{2}}((d-4)h-\sum_{i\geq k+1} e_i)) \to  
H^0(\sH_{2}, \sO_{e_k})
\]
is surjective.
Hence by the exact sequence
\[
0\to \sO_{\sH_{2}}((d-4)h-\sum_{i\geq k} e_i)\to 
\sO_{\sH_{2}}((d-4)h-\sum_{i\geq k+1} e_i)\to \sO_{e_k}\to 0,
\]
we have
$H^1(\sH_{2}, \sO_{\sH_{2}}((d-4)h-\sum_{i=1}^{s} e_i))
\simeq H^1(\sH_{2}, \sO_{\sH_{2}}(d-4)h)$. Since 
it is easy to see that $h^1(\sH_{2},\sO_{\sH_{2}}(d-4)h)=0$,
we have (2).

(3) Since no conic on $A$ intersects all the line on $A$,
$|D_l|$ has no base point.
In case $d=5$, the image of $\Phi_{|D_l|}$ is $\mP^2$ by $(D_l)^2=1$.

Assuming $d\geq 6$, we prove that $D_l$ is very ample.
By (2) and \cite[Theorem 3.1]{DG},
it suffices to prove
that 
\[
h^0(\sH_{2}, \sO_{\sH_{2}}(h-\sum_{j=1}^{d-3} e_{i_j}))=0
\] 
for any set of $d-3$ exceptional curves
$e_{i_1},\dots, e_{i_{d-3}}$.
Assume by contradiction that
there exists an effective divisor 
$L\in |h-\sum_{j=1}^{d-3} e_{i_j}|$ 
for a set of $d-3$ exceptional curves $e_{i_1},\dots, e_{i_{d-3}}$.  
By $\frac{(d-2)(d-3)}{2}-(d-3)\geq 3$, we find at least three $e_i$
such that
$i\not \in \{j_1,\dots,j_{d-3} \}$.
For an $i\not \in \{j_1,\dots,j_{d-3}\}$,
noting $D_l\sim D_{\beta_i}+h-e_{i}¥$, 
$D_l\cdot L=0$, and $L\cdot (h-e_i)>0$, we have
$L\subset D_{\beta_i}$.
This contradicts (\ref{eqnarray:2}) since the number of $i$ 
such that $i\not \in \{j_1,\dots,j_{d-3} \}$ is at least $3$.

We show that
$h^0(\sH_{2}, \sO_{\sH_{2}}(D_l))=d-2$.
By the Riemann-Roch theorem,
$\chi(\sO_{\sH_{2}}(D_l))=d-2$.
Since $h^2(\sH_{2}, \sO_{\sH_{2}}(D_l))=
h^0(\sH_{2}, \sO_{\sH_{2}}(-D_l+K_{\sH_2}))=0$,
we see that
$h^0(\sH_{2}, \sO_{\sH_{2}}(D_l))=d-2$ is equivalent to
$h^1(\sH_{2}, \sO_{\sH_{2}}(D_l))=0$.
Since
$|D_l|$ has no base point, so is $|(d-3)h-\sum_{i\geq k+1} e_i|$.
Thus
the proof that $h^1(\sH_{2}, \sO_{\sH_{2}}(D_l))=0$ 
is almost the same as the above one showing (2)
and we omit it.

(4) follows from
\cite[Proposition 1.1]{Gimi}.
%In case $d\geq 6$, $\sH_2$ is so-called the White surface.
\end{proof}

\begin{rem}
In case of $d=5$,
the morphism defined by $|D_l|$ contracts
three curves $D_{e_i}$ ($i=1, 2, 3$),
which are nothing but the strict transforms of
three lines passing through two of $c_j$.
Namely, the composite $S^2 C\leftarrow \sH_2 \rightarrow \check{\mP}^2$
is the Cremona transformation. 
\end{rem}

\begin{cor}
\label{cor:H_2}
$H^0(\sH_2,\sO_{\sH_2}(i))\simeq
H^0(\check{\mP}^{d-3},\sO_{\check{\mP}^{d-3}}(i))$
for $i=1,2$.
\end{cor}

\begin{proof}
The assertion follows from
Theorem \ref{thm:H_2} (4).
\end{proof} 
The following corollary contains the nontrivial 
result that for a general $[l]\in\sH_{1}$, $D_{l}$ parameterizes
conics which properly intersect $l$.

\begin{cor}\label{noncontenere} 
    For a general $[l]\in\sH_{1}$, $D_{l}$
    does not contain any point corresponding to the line pairs $l\cup m$
    with $[m]\in\sH_{1}$, and hence $D_l$ parameterizes all conics
    which properly interesect $l$.
\end{cor}
\begin{proof} Fix $[m]\in\sH_{1}$ such that $l\cup m$ is a line pair. 
    If $({\overline m},b)$ is the marked line given by $m$
    then we have $d-2$ line pairs $l\cup m$, $l_{1}\cup m$,\ldots
    ,$l_{d-3}\cup m$. Since $L_{b}\sim h$ then $h\cdot D_{l}= d-3$ and
    definitely $[l_{1}\cup m],\ldots
    ,[l_{d-3}\cup m]\in D_{l}$. Thus $[l\cup m]\not\in D_{l}$.
\end{proof}

%&&&&&&&&&&&&&&&CONICS ON A&&&&&&&&&&&&&&&&&&&&&&&&&&&&&&&&&&&&&&&&

%%%%%%%%%%%%%%%%%%%%%%%%%%%%%%%%%%%%%%%%%%%%%%%%%%%%%%%%%%%%%%%%%%%%%%%%%%%%%%%%%%%%%%%%%%%%%%%%%%%%%%%%%%%%%%%%%%%%%%%%%%%%%%%%%%%%%%%%%%%%%%%%%%%%%%%%%%%%%%%%
\section{Varieties of power sums for special quartics $F_4$}
\label{section:VSPA}
~
In Proposition \ref{fuoripiatto} we have seen that 
$\psi\colon \sU_{2}\rightarrow A$ is finite and flat 
outside $\cup_{i=1}^{n}\beta'_{i}$. We can modify 
the morphism $\psi\colon \sU_{2}\rightarrow A$ to obtain a finite one.
See Proposition \ref{finitezzafinale}, which is the goal of 
the subsection \ref{subsection:finite}.

\subsection{Finiteness of 
$\widetilde{\psi}\colon \widetilde{\sU}_2\to \widetilde{A}$}
\label{subsection:finite}~

Let
$\rho\colon{\widetilde{A}}\rightarrow A$ be the blow-up along
$\cup_{i=1}^{s}\beta_{i}'$. We add the following piece of notation:
\begin{nota}
\begin{enumerate}[$(1)$] 
\item
$E_i:=\rho^{-1}(\beta'_{i})$. 
By Lemma \ref{sullebisecanti},
$E_{i}\simeq \mP^{1}\times\mP^{1}$,
%\item
%$f_{i}:=$ a general fiber of $\rho_{\vert E_{i}}\colon
% E_{i}\rightarrow\beta'_{i}$, 
%\item $\gamma_{i}:=$ a general fiber of the other projection
%$E_{i}\rightarrow\mP^{1}$, 
\item $\widetilde{E}_{C}:=$ the strict transform of $E_{C}$, and
\item 
$\widetilde{\zeta}_{ij}:=$ the strict transform of 
the fiber ${\zeta}_{ij}$ of ${E}_C$ over $p_{ij}\in C\cap\beta_{i}$,
\end{enumerate}
where $i=1,\ldots, s$ and $j=1,2$.
\end{nota}

The domain of the finite morphism is ${\widetilde \sU_{2}}:=
\sU_{2}\times_A{\widetilde A}$; in other words,
$\widetilde \sU_{2}$ is the blow-up of $\sU_2$ along 
$\Gamma:=\sU_2\cap (\cup_{i=1}^{s}\beta'_{i}\times \sH_2$). 
We obtain that the natural morphism
${\widetilde\sU_{2}}\rightarrow {\widetilde A}$ is finite after a
local analysis of the morphism ${\sU_{2}}\rightarrow { A}$ in the
neighborhood of $\Gamma$.
 
It is easy to describe $\Gamma$ set-theoretically.
Note that, by Proposition \ref{prop:Cd} (5), 
there are $d-4$ lines $\alpha_{i1},\dots,\alpha_{id-4}$ distinct 
from $\beta_i$ and
intersecting both $C$ and $\beta_i$ outside $C\cap \beta_i$.
Set $t_{ik}:=\alpha_{ik}\cap C$.
Corresponding to $\alpha_{ik}$,
there are two marked conics $(\alpha_{ik}\cup \beta_i; p_{i1}, t_{ik})$
and $(\alpha_{ik}\cup \beta_i; p_{i2}, t_{ik})$.
We denote by $\xi_{ijk}$ the conics on $A$ corresponding to
$(\alpha_{ik}\cup \beta_i; p_{ij}, t_{ik})$, 
where $i=1,\dots, s$, $j=1,2$, and $k=1,\dots, d-4$.
Let $D_{\beta_i}$ be as in the proof of Theorem \ref{thm:H_2}.
Now we can state that 
$\Gamma$ is set-theoretically
the union of
$\beta'_i\times e_i\,$,
\[
\Gamma_i:=
\{(x, [q])\mid [q]\in D_{\beta_i}, x=q\cap \beta'_i\},
\]
which is a section of $\mu$ over $D_{\beta_i}$,
and 
\[
\Gamma_{ijk}:=\beta'_i\times [\xi_{ijk}]\, (i=1,\dots, s,
\, j=1,2,\,
k=1,\dots,d-4).
\]

The conic $\xi_{ijk}$ does not
belong to $e_i$ by the choice of marking.
Moreover, we have the following:

\begin{lem}
\label{lem:0}
The conic $\xi_{ijk}$ does not belong to $D_{\beta_i}$.
\end{lem}

\begin{proof}
We consider the projection of $B$ from a bi-secant line $\beta_i$
(see Proposition \ref{prop:proj1} (1)).
Let $C'\subset Q$ be the image of $C$ by this projection
and $p'_{ij}$ the point of $C'$ corresponding to $p_{ij}$,
where $p_{ij}$ is one of the two point of $C\cap \beta_i$.
By this projection, 
the line $\alpha_{ik}$ maps to a point, which we denote by $s_{ik}$.
Let $F$ be the exceptional divisor of the blowing up along $\beta_i$, 
and $F'$ the image of $F$ on $Q$. 
We say a ruling of $F'\simeq \mP^1\times \mP^1$ is
horizontal if it does not come from a fiber of $F\to \beta_i$. 
Note that the image $q'\subset Q$ 
of a general conic $q$ 
belonging to $D_{\beta_i}$ is a bi-secant line 
of $C'$.
Thus,
if $[\xi_{ijk}]\in D_{\beta_i}$,
then $\xi_{ijk}$ would also correspond to a bi-secant line of $C'$,
which must be the horizontal ruling of $F'$
through $p'_{ij}$ and $s_{ik}$.
By inductive construction of $C$, however,
we can prove that
$p'_{ij}$ and $s_{ik}$ do not lie on a horizontal ruling
(cf. the proof of Lemma \ref{sullebisecanti}).
Thus we have the claim. 
\end{proof}

We can conclude that 
all of $\beta'_i\times e_i$, $\Gamma_i$ and $\Gamma_{ijk}$
are disjoint
$(i=1,\dots, s,
\, j=1,2,\,
k=1,\dots,d-4)$.
%Moreover,
%$\Gamma$ is a local complete intersection curve outside
%$\cup_i \beta'_i\times e_i$
%since so is $\cup_{i=1}^{s}\beta'_{i}\times \sH_2$ in $A\times \sH_2$,
%and $\Gamma$ is one-dimensional outside $\cup_i \beta'_i\times e_i$.

The next proposition contains the final finiteness result we need. 

\begin{prop}\label{finitezzafinale}
    $\widetilde \sU_2$ is Cohen-Macaulay and the natural morphism 
    ${\widetilde\psi}\colon{\widetilde
    \sU_{2}}\rightarrow {\widetilde A}$ is finite $($of degree
    $n:=\frac{(d-1)(d-2)}{2}$$)$.
In particular,
${\widetilde\psi}$ is flat.
\end{prop}

\begin{lem}
\label{lem:Gamma}
$\Gamma$ is a reduced scheme and $\sU_2$ is smooth along $\Gamma$.
\end{lem}

First we finish the proof of Proposition \ref{finitezzafinale}
by admitting this lemma:

\begin{proof}[Proof of Proposition $\ref{finitezzafinale}$]
By Lemma \ref{lem:Gamma},
the morphism $\widetilde{\sU}_2\to \sU_2$
is the blow-up along the reduced subscheme $\Gamma$ contained in
the smooth locus of $\sU_2$.
The subscheme $\beta'_i\times e_i$ is a Cartier 
divisor of $\sU_2$,
thus $\widetilde{\sU}_2\to \sU_2$ is isomorphic over 
$\beta'_i\times e_i$.
The curve $\Gamma_{ijk}$ is smooth and
the curve $\Gamma_i$ has only planar singularities since so is $D_{\beta_i}$.
Thus $\widetilde{\sU}_2$ is Cohen-Macaulay
since so is $\sU_2$. 

We have only to prove that
$\widetilde \psi$ is finite.
By Proposition \ref{fuoripiatto},
$\widetilde \psi$ is finite outside $\cup_{i} E_i$.
Note that
$\widetilde{\psi}^{-1} (E_i)$ is nothing but
the inverse images of $\beta'_i\times e_i$, $\Gamma_i$ and $\Gamma_{ijk}$
by $\widetilde{\sU}_2\to \sU_2$,
all of which are $\mP^1$-bundles over curves
and are mapped to $E_i$ finitely.
Hence we are done.
\end{proof}

\begin{proof}[Proof of Lemma $\ref{lem:Gamma}$]
We study $\sU_2$ locally along $\Gamma$.

Let $q$ be a conic on $A$ belonging to $D_{\beta_i}$.
Then,
by Proposition \ref{prop:Cd1} (5),
Lemma \ref{lem:0} and the fact that $D_{\beta_i}\cap e_i=\emptyset$
(see the proof of Theorem \ref{thm:H_2} (\ref{eqnarray:1})),
we see that $q$ is smooth near $\beta'_i$
and intersects $\beta'_i$ transversely.
This implies that $\widetilde{\sU}_2$ is smooth along $\Gamma_i$.
Note that, near $\Gamma_i$, 
the morphism $\psi\colon \sU_2\to A$ is finite, hence flat. 
Since $\Gamma$ is the pull-back of $\beta'_i$ near $\Gamma_i$ and
$\Gamma_i$ is not contained in the ramification locus of $\psi$,
it holds that
$\Gamma$ is reduced along $\Gamma_i$.

Let $q$ be the fiber of 
$\sU_2\to \sH_2$ over $[\xi_{ijk}]$ or 
a point of $e_i$.
Note that $q$ is a conic on $A$
and has only nodes as its singularities.
We show that
$h^1(\sN_{q/A})=0$ and 
the natural map $H^0(\sN_{q/A})\to H^0(T^1_p)\simeq \mC$
is surjective, where $p$ is any node of $q$
and $T^1_p$ is the local deformation space of $p$.
As in the proof of \cite[Proposition 1.1]{HH},
this implies that $\sH_2$ coicides with
the Hilbert scheme of conics on $A$
at $[\xi_{ijk}]$ or 
a point of $e_i$,
and $\sU_2$ is smooth near $q$.

First we treat the case where $q=\xi_{ijk}=\alpha'_{ik}\cup \beta'_i\cup
\zeta_{i,3-j}$.
Note that 
$\sN_{\alpha'_{ik}/A}\simeq \sO_{\mP_1}\oplus \sO_{\mP_1}(-1)$,
$\sN_{\beta'_i/A}\simeq \sO_{\mP_1}(-1)^{\oplus 2}$,
and
$\sN_{\zeta_{i,3-j/A}}\simeq \sO_{\mP_1}\oplus \sO_{\mP_1}(-1)$.
We apply \cite[Theorem 4.1]{HH} by setting
$X=\xi_{ijk}$, $C=\beta'_i$ and $D=\alpha'_{ik}\cup \zeta_{i,3-j}$.
% and $S=(\alpha'_{ik}\cap \beta'_i)\cup (\zeta_{i,3-j}\cap \beta'_i)$.
We check the conditions a) and b) of [ibid.].
The condition a) clearly holds.
The condition b) follows from the following two facts:
\begin{enumerate}[(1)]
\item
let $F$ be the exceptional divisor of 
the blow up of $B$ along $\alpha_{ik}$.
Note that $F\simeq \mP^1\times \mP^1$.
We call a fiber of $F\to \mP^1$ in the other direction to
$F\to \alpha_{ik}$ a horizontal fiber.
Then the intersection points of 
the strict transform of $C$ and $F$, and 
the strict transform of $\beta_i$ and $F$
do not lie on a common horizontal fiber.

This can be proved by the inductive construction of $C=C_d$
in a similar fashion to the proof of Lemma \ref{sullebisecanti}, 
or by a straightforward dimensional computation as the one of 
Proposition \ref{prop:Cd1} (2), and
\item
let $G$ be the exceptional divisor of 
the blow up of $A$ along $\zeta_{i,3-j}$.
Note that $G\simeq \mF_1$.
Then the intersection points of 
the strict transform of $\beta'_i$ and $G$
does not lie on the negative section of $G$.

Indeed, since $E_C\cdot \zeta_{i,3-j}=-1$,
the intersection of $G$ and the strict transform of $E_C$
is the negative section of $G$.
On the other hand, the strict transforms of 
$E_C$ and $\beta'_i$ are disjoint.
\end{enumerate}
Thus, by \cite[Theorem 4.1]{HH}, $\xi_{ijk}$ 
satisfies the desired properties.

Secondly, we treat the case $q$ is a fiber over a point of $e_i$.
Note that $\overline{q}=\beta_i\cup \alpha$,
where $\alpha$ is a line intersecting $\beta_i$.
Denote by $\alpha'$ the strict transform of $\alpha$.
We make the following case division:
\begin{enumerate}[(a)]
\item
$\alpha\cap C=\emptyset$ and $\sN_{\alpha/B}=\sO_{\mP^1}^{\oplus 2}$.
\item
$\alpha\cap C=\emptyset$ and $\sN_{\alpha/B}=\sO_{\mP^1}(-1)\oplus
\sO_{\mP^1}(1)$.
\item
$\alpha=\alpha_{ik}$ for some $k$.
\item
$\alpha$ passes through a point of $\beta_i\cap C$. 
\end{enumerate}

In the case (a) or (b),
it is easy to see 
the proof of \cite[Theorem 4.1]{HH} works as above
by setting $X=q$, $C=\beta'_i$ and $D=\alpha'$.
In the case (c) or (d), we need to modify the proof of [ibid.].
We only treat the case (c) since we can treat the case (d) similarly.
Note that $q=\beta'_i\cup \alpha'_{ik}\cup
\gamma_{ik}$, where $\gamma_{ik}$ is the fiber of $E_C$ over
$t_{ik}$ Note that $C$ is smooth.
%As above (by inductive construction of $C$),
By \cite[Corollary 3.2]{HH} and simple dimension count,
we can describe the restrictions of the normal bundle
$\sN_{q/A}$ to the components of $q$ as follows:
\[
\sN_{q/A|\beta'_i}=\sO_{\mP^1}
\oplus \sO_{\mP^1}(-1),
\sN_{q/A|\alpha'_{ik}}=\sO_{\mP^1}
\oplus \sO_{\mP^1}(1),\
\text{and} \
\sN_{q/A|\gamma{ik}}=\sO_{\mP^1}^{\oplus 2}.
\]
Set $C=\beta'_i\cup \gamma_{ik}$ and $D=\alpha'_{ik}$.
As in \cite[Theorem 4.1]{HH}, set $S:=C\cap D$.
By the description of 
$\sN_{q/A|\beta'_i}$, 
$\sN_{q/A|\alpha'_{ik}}$,
and 
$\sN_{q/A|\gamma{ik}}$,
it holds that
$H^1(\sN_{q/A|C})=H^1(\sN_{q/A|D})=\{0\}$.
Moreover, considering the tautological linear systems of
$\mP(\sN_{q/A|\beta'_i})$, 
$\mP(\sN_{q/A|\alpha'_{ik}})$,
$\mP(\sN_{q/A|\gamma{ik}})$,
and $\mP(\sN_{q/A})$,
we see that 
$H^0(\sN_{q/A|C})\oplus H^0(\sN_{q/A|D})\to H^0(\sN_{q/A|S})$
is surjective.
Thus $h^1(\sN_{q/A})=0$ holds.
By \cite[Corollary 3.2]{HH} again, we have the following exact sequences
(cf. \cite[(3) in the proof of Theorem 4.1]{HH}):
\[
0\to \sO_{\mP^1}(-1)\oplus \sO_{\mP^1}(-2) 
\to \sN_{q/A|\beta'_i}\to \sN_{q/A|S}\to 0,
\]
\[
0\to  \sO_{\mP^1}(-1)\oplus \sO_{\mP^1}(-2)
\to \sN_{q/A|\alpha'_{ik}}
\to \sN_{q/A|S}\to 0,
\]
\[
0\to  \sO_{\mP^1}^{\oplus 2}(-1)\to \sN_{q/A|\gamma{ik}}\to \sN_{q/A|S}\to 0.
\]
Thus we can consider that
$H^0(\sN_{q/A|C})$ and $H^0(\sN_{q/A|D})$
is a subspace of $H^0(\sN_{q/A|S})$.
By \cite[(2) in the proof of Theorem 4.1]{HH},
we see that
$H^0(\sN_{q/A|C})\to H^0(T^1_p)$ and
$H^0(\sN_{q/A|D})\to H^0(T^1_p)$ are surjective.
Moreover,
considering the tautological linear systems of
$\mP(\sN_{q/A|\beta'_i})$, 
$\mP(\sN_{q/A|\alpha'_{ik}})$,
$\mP(\sN_{q/A|\gamma{ik}})$,
and $\mP(\sN_{q/A})$,
we see that
the kernels of $H^0(\sN_{q/A|C})\to H^0(T^1_p)$ and
$H^0(\sN_{q/A|D})\to H^0(T^1_p)$ does not coincide
for any $p\in S$. 
Thus any non-zero element of $H^0(T^1_p)\simeq \mC$
comes from that of $H^0(\sN_{q/A|C})\cap H^0(\sN_{q/A|D})$
as in the end of the proof of 
\cite[Theorem 4.1]{HH}. This implies that 
the natural map $H^0(\sN_{q/A})\to H^0(T^1_p)$
is surjective for any $p\in S$.

Note that, near $e_i$,
the family $\sU_2\to \sH_2$ is locally a deformation of a node
with smooth discriminant locus $e_i$.
Thus a local computation shows that
$\Gamma$ is reduced along $\beta'_i\times e_i$.

Now we prove that $\Gamma$ is reduced along $\Gamma_{ijk}$.
We have only to prove that
${\sU}_2\to A$ is unramified along $\Gamma_{ijk}$
since then $\Gamma$ is the \'etale pull-back of $\beta'_i$ near
$\Gamma_{ijk}$, hence is reduced.

Recall that we set
$S=(\alpha'_k\cap \beta'_i)\cup (\zeta_{i,3-j}\cap \beta'_i)$.
By simple dimension count and
the following exact sequence:
\[
0\to \sN_{\beta'_i/A}\to \sN_{\xi_{ijk}/A|\beta'_i}\to T^1_S\to 0,
\]
we can prove that $\sN_{\xi_{ijk}/A|\beta'_i}\simeq \sO_{\mP^1}^{\oplus 2}$.
Thus $H^0(\sN_{\xi_{ijk}/A})\otimes \sO_{\xi_{ijk}}\to \sN_{\xi_{ijk}/A}$
is surjective at a point of $\Gamma_{ijk}$ since
it factor through the surjection
$H^0(\sN_{\xi_{ijk}/A|\beta'_i})\otimes \sO_{\beta'_i}\to 
\sN_{\xi_{ijk}/A|\beta'_i}$.
Thus 
${\sU}_2\to A$ is unramified along $\Gamma_{ijk}$.
\end{proof}

From now on
we assume that $d\geq 6$ 
and we consider $\sH_2\subset \check{\mP}^{d-3}$.

Consider  the following
diagram:

\begin{equation}
\label{iniezionehilbert}
\xymatrix{ & \widetilde{\sU_{2}} \ar[dl]_{{\widetilde{\mu}}}
\ar[dr]^{{\widetilde\psi}}\\
 \sH_{2} &  & {\widetilde A}. }
\end{equation}

\begin{defn}\label{attached}
    Let $\widetilde a$ be a point of $\widetilde A$.
    We say that $[\widetilde \psi^{-1}(\widetilde a)]\in 
\Hilb^{n} \check{\mP}^{d-3}$ is {\it{the cluster of conics attached to}} $\widetilde a$ and denote it by $[\sZ_{\widetilde a}]$. A conic $q$ such that $[q]\in \Supp \sZ_{\widetilde a}$ is 
called {\it{a conic attached to $\widetilde a$.}}
 \end{defn}

\begin{rem}
Though we do not need it later,
we describe the fiber of $\widetilde \psi$ over a general point 
${\widetilde a}\in E_i$ for some $i$
for reader's convenience.
In other words, we exhibit $n$ conics 
attached to $\widetilde a$.

Set $a:=\rho({\widetilde a})\in A$ and
    $b:= f(a)\in \beta_i$. We use notations of Proposition \ref{fuoripiatto}. 
    Since $\deg C_b=d-2$,
the number of bi-secant conics through $b$ not belonging to
the family $e_i$ is given by the number of double points of $C_b$, 
which is $\frac{(d-3)(d-4)}{2}$. Moreover
$2(d-4)$ conics $\xi_{ijk}$ through $a$.

The number of remaining conics is $3=n-\frac{(d-3)(d-4)}{2}-2(d-4)$. 
Such conics will belong to $e_i$.
We look for three such conics.
Let $S_i$ be the strict transform on $\widetilde{A}$ of 
the locus of lines intersecting $\beta_i$. 
Then it is easy to see that $S_{i|E_i}$ does not contain 
any fiber $\gamma_{i}$ of the second projection 
$\sigma_{i}\colon E_i\rightarrow\mP^{1}$.
Moreover $S_{i|E_i}\sim2\gamma_i+3f_{i}$,
where $f_{i}$ is a fiber of $E_i\to \beta'_i$.
Let $\gamma'_{i}$ be the fiber of $\sigma_{i}$ through $\widetilde{a}$.
Then $\gamma'_{i}¥$ intersect $S_i$ at three points.
Corresponding to these three points, there are
three lines on $B$ intersecting $\beta_i$.
Denote by $l_1$, $l_2$ and $l_3\subset A$ the strict transforms of these three 
lines. Then $\beta'_i\cup l_j$ $(j=1,2,3)$ are the conics on $A$ what we want.
\end{rem}

By Proposition \ref{finitezzafinale} and 
the universal property of Hilbert schemes, we obtain
a naturally defined map $\Psi\colon
{\widetilde A}\rightarrow{\rm{Hilb}}^{n}\check{\mP}^{d-3}$.
This is clearly injective because 
$n$ conics attached to a point ${\widetilde a}\in
    {\widetilde A}$ uniquely determines $\widetilde a$.

To understand the image of $\Psi$,
we construct the special quartic hypersurface
which live in the dual projective space to the ambient of $\sH_2$.

%%%%%%%%%%%PIATTEZZA:FINE%%%%%%%%%%%%%%%%%%%%%%%%%%%%%%%

%%%%%%%%%%%%%%%%%%%%%%%%%%%%%%%%%%%%%%%%%%%%%%%%%%%%%%%%%

\subsection{Intersection of conics and conics on $A$}
\label{subsection:conics}~

To construct the special quartic hypersurface,
we need the incidence variety defined
by the intersections of conics.

Similarly to (\ref{eq:familyB1}),
we consider the following diagram:
\begin{equation}
\label{eq:family2}
\xymatrix{\sU'_2\subset {\sU}_2\times \sH_2 \ar[d] 
\ar[r]^{(\psi, \mathrm{id})} & A \times \sH_2 \supset \sU_2 \ar[d]\\
%\overline {\sU}_2\times \sH_2 \ar[d]_{({\overline{\mu}},\mathrm{id})} 
%\ar[r]^{(\overline \psi, \mathrm{id})} & B \times \sH_2 \supset
%\overline \sU_2 \ar[d]\\
\widehat{\sD}_2 \subset \sH_2\times \sH_2 \ar[r] & \sH_2,} 
\end{equation}
where $\sU'_2\subset \sU_2\times \sH_2$ is the base change of $\sU_2$
and $\widehat{\sD}_2$ the image of $\sU'_2$ on $\sH_{2}\times \sH_{2}$.
Similarly to the investigation of the diagram (\ref{eq:familyB1}),
we see that the image $\sF'$ in $\sH_2\times \sH_2$ of the inverse image of 
$\cup_{i=1}^{n}\beta'_i\times \sH_2$ is not divisorial nor does not dominate $\sH_2$.
Moreover, 
any component of $\widehat{\sD}_2$ outside $\sF'$ dominates $\sH_2$, and
is divisorial or possibly the diagonal of $\sH_2\times \sH_2$.
Note that dislike the diagram (\ref{eq:familyB1}),
there is no other non-divisorial component in this case.
Compare the proof of Corollary \ref{divisoriali}.
Here we leave the possibility that 
the diagonal of $\sH_2\times \sH_2$ is contained in
the divisorial component of $\widehat{\sD}_2$.
We, however, prove this is not the case in 
Lemma \ref{noselfintersection}. 
%Moreover, 
%except possibly the image of such components and the diagonal of $\sH_2\times \%sH_2$, any other component of $\widehat{\sD}_2$ is a prime divisor. 
%Moreover, no divisorial component of $\widehat{\sD}_2$ dominates $\sH_2$ 
%and any divisorial component of $\widehat{\sD}_2$ dominates $\sH_2$.

Let $\sD_2 \subset \sH_2\times \sH_2$ be the union of the
divisorial components of $\widehat{\sD}_2$ with reduced structure.
$\sD_2$ is Cartier since $\sH_2\times \sH_2$ is smooth.
$\sD_2 \to \sH_2$ is flat 
since $\sD_2$ is Cohen-Macaulay, $\sH_2$ is smooth and 
$\sD_2 \to \sH_2$ is equi-dimensional.  
Let $D_q$ be the fiber of 
$\sD_2\to \sH_2$ over $[q]\in \sH_2$ via the projection to
the second factor.

\subsection{Construction of the special quartics}
~

\begin{lem}
\label{lem:Dq}
$D_q\sim
2(d-3)h-2\sum_{i=1}^{e} e_i$
for a conic $q$.
$D_q$ is a quadric section of $\sH_2\subset \check{\mP}^{d-3}$.
\end{lem}

\begin{proof}
The proof of the former half is almost identical to the one of 
Theorem \ref{thm:H_2} (1).
The latter half follows from Corollary \ref{cor:H_2}.
\end{proof} 

Now we proceed to construct the quartic hypersurface. 
%The construction has many overlaps with that of the Scorza quartic
%as in \cite[\S 9]{DK}.

From now on, we write $\mP^{d-3}=\mP_* V$, where $V$ is 
the $d-2$-dimensional vector space.
The crucial point in the following considerations is the equality:
\begin{equation}
\label{eq:n}
n=\dim S^{2}V.
\end{equation}

By the seesaw theorem,
it holds that $\sD_2\sim p_1^*D_q+p_2^*D_q$.
Consider 
the morphism $\sH_2\times \sH_2$ into 
$\check{\mP}^{d-2}\times \check{\mP}^{d-3}$
defined by $|p_1^*D_l+p_2^*D_l|$, which 
is an embedding since $d\geq 6$.
By Corollary \ref{cor:H_2}, it holds
\[
H^0(\sH_2\times \sH_2,\sD_2)
\simeq H^0(\check{\mP}^{d-3}\times \check{\mP}^{d-3},
\sO(2,2)).
\] 
Therefore $\sD_2$ is the restriction of a unique $(2,2)$-divisor on 
$\check{\mP}^{d-3}\times \check{\mP}^{d-3}$, 
which we denote by $\{\widetilde{\sD}_2=0\}$.
Since $\{\widetilde{\sD}_2=0\}$
is symmetric, 
we may assume the equation $\widetilde{\sD}_2$ is also symmetric.
Actually, the desired quartic is obtained by
restricting $\widetilde{\sD}_2$ to the diagonal and
taking the dual in the sense of Dolgachev (see the appendix), 
but we need more argument for the proof of the main theorem.

For $[q]\in \sH_2$, we denote 
by
$\widetilde{D}_{q}$
the restriction of 
$\widetilde{\sD}_2$ to the fiber over $[q]$.
Note that 
$\widetilde{\sD}_2\in S^2 {V} \otimes S^2 {V}$,
so $\widetilde{\sD}_2$ defines a linear map
$\lambda\colon S^2 \check{V} \simeq (S^2 {V})\check{\empty}\to S^2 V$.
Let $H_q$ be a linear form on $\check{V}$
corresponding to $q$.
It holds that
$\lambda(H^2_q)=\widetilde{D}_q$ up to scalar,
so we may choose $H_q$ such that
$\lambda(H^2_q)=\widetilde{D}_q$ holds.  
We prove that $\lambda$ is an isomorphism.

\begin{lem}\label{noselfintersection}
$\sD_2$ does not contain the diagonal of $\sH_2\times \sH_2$.
In particular we have the following:

let $\widetilde a$  be a general point of $\widetilde{A}$ and
$q_{1}, q_{2},\ldots ,q_{n} \in\sH_{2}$
the conics attached to $\widetilde a$. Then 
$$D_{q_i}([q_i])\neq 0$$
\noindent
for $1\leq i\leq n$.
\end{lem}
\begin{proof}
Here we assume $d\geq 3$. 
It suffices to prove that
$D_q([q])\neq 0$
for a general $[q]\in \sH_2$.
This is equivalent to that
the image ${D}^{\flat}_q$ on $\overline{\sH}_2$ of $D_q$ does not contain
$[\overline{q}]$. Note that $D^{\flat}_q$ is the closure
of the locus of multi-secant conics of $C$ intersecting 
properly $\overline{q}$.
Now the assertion follows from the inductive construction of $C_d$ from
$C_{d-1}\cup \overline{l}$.
From now on, we denote ${D}^{\flat}_{q}$ for $C_d$ by ${D}^{\flat}_{q,d}$.
If $d=3$, then $D_q\sim 0$, thus 
the assertion trivially true.
If ${D}^{\flat}_{q',d-1}([\overline{q}'])\neq 0$ 
for a general multi-secant conic $\overline{q}'$ of $C_{d-1}$,
then ${D}^{\flat}_{q,d}([\overline{q}])\neq 0$ for 
a general multi-secant conic $\overline{q}$ of $C_{d}$. 
\end{proof}
 Let $\widetilde{a}$ be a general point of $\widetilde{A}$ and 
    $q_{1},\ldots, q_{n}$ are the conics attached to $\widetilde{a}$.
    By the definition of $\widetilde{D}_{q_i}$ and 
    generality of $\widetilde{a}$, 
we have the following:
\begin{equation}
\label{eq:1'}
\text{$\widetilde{D}_{q_j}([q_i])=0$ $(j\not =i)$
and $\widetilde{D}_{q_i}([q_i])\neq 0$}.
\end{equation}
 
    (\ref{eq:1'}) implies 
    $\widetilde{D}_{q_1},\dots, \widetilde{D}_{q_n}$ 
    are linearly independent, and, by 
    (\ref{eq:n}), they span the vector space $S^2 V$.
    Thus $\lambda$ is an isomorphism.

The inverse $\lambda^{-1}\colon S^2 V \to S^2 \check{V}$
defines an element $\check{\sD}_2\in S^2 \check{V}\otimes
S^2 \check{V}$.
We consider 
the polarization map $\mathrm{pl}_2\colon S^2 \check{V}\to \mathrm{Sym}_2 V$
(see the appendix).
We show that 
$\widetilde{U}:=\mathrm{pl}_2\otimes \mathrm{pl}_2 (\check{\sD}_2)
\in \mathrm{Sym}_2 V\otimes \mathrm{Sym}_2 V\subset {\check{V}}^{\otimes 4}$
is actually contained in $\mathrm{Sym}_4 V$.
This will implies that
$\mathrm{pl}_2\otimes \mathrm{pl}_2 (\widetilde{\sD}_2)$
is the image of a quartic form $\in S^4 \check{V}$ by $\mathrm{pl}_4$. 
 
The following argument is almost identical with
the proof of \cite[Theorem 9.3.1]{DK}
(The identification will be clearer by constructing
the theta characteristic on $\sH_1$ in the forthcoming paper).
Let $l$ be a general line on $A$ and
$l_1,\dots, l_{d-2}$ the lines intersecting $l$.
Note that
$l_1,\dots, l_{d-2}$ correspond to lines on $B$ intersecting
both $C$ and the image $\overline{l}$ 
of $l$ on $B$ except those through $C\cap \overline{l}$.
Thus the number of such lines is $d-2$.
Since $l$ is general, 
so are $l_1,\dots, l_{d-2}$.
We have $d-2$ reducible conics $r_1:=l\cup l_1,\dots, r_{d-2}:=l\cup l_{d-2}$.
It holds that $D_{r_i}=D_{l}+D_{l_i}$.
By Corollary \ref{cor:H_2},
$\widetilde{D}_l$ and $\widetilde{D}_{l_i}$ are defined 
by linear forms $L$ and
$L_i$.
We may assume that
$\lambda(H_{r_i}^2)=\widetilde{D}_{r_i}=L_{i}L$.
By Corollary \ref{noncontenere}, 
$L_i([r_i])\not =0$ and
$L_i([r_j])=0$ for $i\not =j$.
In other words,
it holds
$\langle L_i, H_i\rangle\not =0$ and
$\langle L_i, H_j\rangle =0$ for $i\not =j$,
where $\langle \empty, \empty\rangle$ is the
natural dual pairing.
Thus $L_1,\dots, L_{d-2}$ and 
$H_{r_1}, \dots, H_{r_{d-2}}$ span
$\check{V}$ and $V$, respectively
since $\dim \check{V}=d-2$.
Moreover, $\{H_{r_i}\}$ and $\{L_i\}$ are dual to each other.
Choose coordinates of $V$ and $\check{V}$ such that
$H_{r_i}$ and $L_i$ are coodinate hyperplanes $\{x_i=0\}$
and $\{u_i=0\}$ respectively.
Set $L=\sum a_i u_i$.
For any $y=(y_1,\dots,y_{d-2})\in V$,
we have
$\lambda(\sum y_i x_i^2)=(\sum a_i u_i)(\sum y_i u_i)$
by $\lambda(H^2_{r_i})=L_iL$.
Considering 
$\widetilde{U}
\in \check{V}^{\otimes 4}$,
this implies that
$\widetilde{U}(L,y, x,x)=\sum y_i x_i^2=P_y (\sum x_i^3)$,
where $x=(x_1,x_2,\dots,x_{d-2})$ and $P_y$ is the polar with
respect to $y$ (see the appendix).
Thus we have $\widetilde{U}(L,y, x,z)=\sum y_i x_i z_i$
for $z=(z_1,z_2,\dots,z_{d-2})$, hence
$\widetilde{U}(L,y, x,z)$ is symmetric for $y$, $x$ and $z$.
Since 
$\widetilde{U}\in \mathrm{Sym}_2 \check{V}\otimes \mathrm{Sym}_2 \check{V}$
and $\widetilde{\sD}_2$ is symmetric,
we have shown that $\widetilde{U}\in \mathrm{Sym}_4 \check{V}$.

Let $F_4$ be the quartic form associated to $\widetilde{U}$,
namely, $F_4:=\widetilde{U}(x,x,x,x)$.
By the construction, it holds
\begin{equation}
\label{eq:F_4}
P_{\widetilde{D}_q}(F_4)=H_q^2.
\end{equation}

By the theory of polarity (see the appendix),
we can interpret what we have done as follows:

$\lambda^{-1}=\mathrm{ap}^2_{F_4}$. 
Since $\lambda^{-1}$ is an isomprphism,
$F_4$ is non-degenerate.

\subsection{Proof of the main theorem}
\label{subsection:main}
\begin{thm}\label{diretto}
$\Ima \Phi$ is an irreducible component of
\[
\VSP(F_4, n;\sH_2):=
\overline{
\{([H_1],\dots, [H_n])\mid [H_i] \in \sH_2\}}
\subset \VSP(F_4,n).
\]

\end{thm}

\begin{proof}
Set 
\[
Z:=\{([H_{1}],\ldots, [H_{n}])
    \in \Hilb^{n} \check{\mP}^{d-3} \mid H^{4}_{1}+\ldots
 +H^{4}_{n}=F_4, [H_i]\in \sH_2\}.
\]

    %First I prove $\mathrm{Image}\, \Psi\subset \overline{Z}$.
    %and $\rank F_2=d$. 
    For a general point $\widetilde{a}$ and conics  
    $q_{1},\ldots, q_{n}$ attached to $\widetilde{a}$, 
we have (\ref{eq:1'}).
Conversely, $n$ conics $q_i$ satisfying (\ref{eq:1'}) and
the assumptions (1)--(3) of Lemma \ref{lem:GeomA}
determine a point of $\widetilde{A}$.
Note that the assumptions (1)--(3) of Lemma \ref{lem:GeomA}
are open conditions.
Thus   
we have only to prove that (\ref{eq:1'}) is equivalent to
\begin{equation}
\label{eq:2'}
    \text{$\alpha_1 H^{4}_{q_1}+\ldots +\alpha_n H^{4}_{q_n}=F_4$
    with some nonzero $\alpha_i \in \mC$.} 
\end{equation}

We see that (\ref{eq:2'}) is equivalent to
\begin{equation}
\label{eq:3'}
\text{If $\{G=0\}\subset \check{\mP}^{d-3}$ is any quartic
through
    $[q_1], \cdots, [q_n]$, then $P_{F_4}(G)=0$.}
\end{equation}

Indeed, 
by the apolarity pairing,
$\langle G, H^4_{q_i} \rangle=0\Leftrightarrow  G([q_i])=0$,
thus, the assumption on $G$ is equivalent to
$G\in \langle H^4_{q_1},\dots,H^4_{q_n}\rangle^{\perp}.$
Therefore
(\ref{eq:2'}) is equivalent to
$\langle H^4_{q_1},\dots,H^4_{q_n}\rangle^{\perp}\subset 
\langle F_4 \rangle^{\perp}.$
Since $F_4$ is non-degenerate,
this is equivalent to (\ref{eq:2'}).    

We show (\ref{eq:1'}) implies (\ref{eq:3'}).
If (\ref{eq:1'}) holds, then
    $\widetilde{D}_{q_i}$ ($i\not =1$)
    generate the space of quadric forms passing through $[q_1]$,
    we may write 
$G=Q_2\widetilde{D}_{q_2}+\cdots+Q_n\widetilde{D}_{q_n}$,
    where $Q_i$ are quadratic forms on $\check{\mP}^{d-3}$.
    By $G([q_i])=0$ for $i\not =1$,
    we have $Q_i([q_i]) \widetilde{D}_{q_i}([q_i])=0$.
    $\widetilde{D}_{q_i}([q_i])\not=0$ implies that $Q_i([q_i])=0$. Thus
    $Q_i$ is a linear combination of
    $\widetilde{D}_{q_j}$ $(j\not =i)$.
    Consequently,
    $G$ is a linear combination of
    $\widetilde{D}_{q_i}\widetilde{D}_{q_j}$ $(1\leq i<j\leq n)$.
    Thus $P_{F_4}(G)=0$ follows from 
    that
   \[
   P_{F_4}(\widetilde{D}_{q_i}\widetilde{D}_{q_j})=
   P_{H_{q_i}}(\widetilde{D}_{q_j})=\widetilde{D}_{q_j}([q_i])=0.
   \]
Finally we show (\ref{eq:2'}) implies (\ref{eq:1'}).
By (\ref{eq:2'}), it holds
\[
H_{q_i}^2=P_{\widetilde{D}_{q_i}}(F_4)=
\sum \alpha_j \langle \widetilde{D}_{q_i}, H_{q_j}^4 \rangle H_{q_j}^2.
\]
Since $\widetilde{D}_{q_i}$ are linearly independent,
so are $H_{q_j}^2$. 
Thus (\ref{eq:1'}) holds.
\end{proof}

\begin{defn}
We say $\Ima \Phi$ is the {\em{main component}} of $\VSP(n,F_4;\sH_2)$.
\end{defn}

The following proposition characterizes the main component of 
$\VSP(n,F_4;\sH_2)$:
%which will play a crucial role in \ref{subsection:moduli}:

\begin{prop}
\label{lem:main}
Let $(\sH_2^k)^o$ and 
$(\Hilb^k \check{\mP}^{d-3})^o$ $(k\in \mN)$ be the complements of
all the small diagonals of $\sH_2^k$ $(k$ times product of $\sH_2)$ and $\Hilb^k \check{\mP}^{d-3}$
respectively.
Set
\[
\VSP^o(F_4,n;\sH_2):=
\{([H_1],\dots, [H_n])\mid [H_i] \in \sH_2, H_1^m+\cdots+H_n^m=F_4\}.
\]
Let $V^o$ be the inverse image of
$\VSP^o(F_4,n;\sH_2)$ by the natural projection
$(\sH_2^n)^o\to (\Hilb^n \check{\mP}^{d-3})^o$.
Let $(\sH_2^n)^o\to (\sH_2^2)^o$ be the projection
to any of two factors.
Then a component of $V^o$ dominating $\sD_2$
dominates the main component of 
$\VSP(F_4,n;\sH_2)$.
In particular, the main component is uniquely identified
from $\sD_2$.
\end{prop}

\begin{proof}
Let $([q_1],[q_2])\in \sD_2\cap (\sH_2^2)^o$ be a general point
and $\{q_i\}$ ($i=1,\dots,n$) any set of mutually conjugate
$n$ conics including $q_1$ and $q_2$.
Since $q_1$ and $q_2$ are general, we may assume that
all the $q_i$ are general. 
By Lemma \ref{lem:GeomA} and Theorem \ref{diretto},
it suffices to prove that
$q_1,\dots,q_n$ satisfies 
the conditions (1)--(3) of Lemma \ref{lem:GeomA}.

(1). Let $\overline{r}_1$ and $\overline{r}_2$ are mutually intersecting
smooth conics on $B$
and $\overline{r}_3$ a line pair on $B$ intersecting
both $\overline{r}_1$ and $\overline{r}_2$.
Since the Hilbert scheme of conics on $B$ is $4$-dimensional,
the pair of $\overline{r}_1$ and $\overline{r}_2$ depends on $7$ parameters.
If we fix $\overline{r}_1$ and $\overline{r}_2$,
then $\overline{r}_3$ depends on $1$ parameter. 
Thus the configuration $\overline{r}_1$, $\overline{r}_2$, $\overline{r}_3$
depends on $8$ parameters.
Fix $\overline{r}_1$, $\overline{r}_2$ and $\overline{r}_3$.
We count the number of parameters of $C_d$ such that
$C_d$ intersects each of $\overline{r}_i$ ($i=1,2,3$) twice.
The number of parameters is 
$h^0((\sO_{\mP^1}(d-1)\oplus \sO_{\mP^1}(d-1))\otimes \sO_{\mP^1}(-6))+6=2d-12+6=2d-6$,
where $+6$ means the sum of the numbers of parameters of two points on 
$\overline{r}_i$ ($i=1,2,3$).
By $2d-6+8=2d+2$,
a general $C_d$ has $2$-dimensional pairs of mutually intersecting bi-secant 
conics which intersect at least one bi-secant line pair of $C_d$.
Thus general pairs of mutually intersecting bi-secant conics of $C_d$, 
which form a $3$-dimensional family, do not intersect
a bi-secant line pair of $C_d$.
  
(2).
Assume by contradiction that  
$\overline{q}_i$, $\overline{q}_j$ and $\overline{q}_k$
pass through a point $b$, and
$\overline{q}_l$ does not pass through $b$ but intersects a line
through $b$. 
Then by the double projection from $b$,
$\overline{q}_l$ is mapped to a line through 
the three singular points of the image of $C_b$ corresponding to  
$\overline{q}_i$, $\overline{q}_j$ and $\overline{q}_k$.
Thus we have only to prove that for a general point of $b$ on $B$,
three double points of the image of $C_b$ do not lie on a line.

Fix a general point $b\in B$.
Let $\overline{r}_1$, $\overline{r}_2$, $\overline{r}_3$ be three conics on $B$ through $b$
such that by the double projection from $b$,
they are mapped to three colinear points on $\mP^2$.
%Assume that $d\geq 6$ (if $d=5$, then direct verification).
The number of parameters of $C_d$'s intersecting each of $\overline{r}_i$ twice
is 
$h^0((\sO_{\mP^1}(d-1)\oplus \sO_{\mP^1}(d-1))\otimes \sO_{\mP^1}(-6))=2d-12$
since $h^1((\sO_{\mP^1}(d-1)\oplus \sO_{\mP^1}(d-1))\otimes \sO_{\mP^1}(-6))=0$
Note that the number of parameters of
$\overline{r}_1$, $\overline{r}_2$, $\overline{r}_3$
is $5$ since
that of lines in $\mP^2$ is $2$,
and that of three points on a line
is $3$.
Thus the number of parameters of $C_d$'s 
such that its image of the double projection from $b$
has three colinear double points is at most $2d-1$.
Hence a general $C_d$ does not satisfy this property.

(3).
Let $r_1$ and $r_2$ be a general pair of mutually conjugate conics on $A$
such that
$\overline{r}_1$ and $\overline{r}_2$ are smooth, and
$\overline{r}_1$ and $\overline{r}_2$ 
intersect at a point on $C\cup \cup_i \beta_i$.
Such general pairs of conics $r_1$ and $r_2$ 
form a two-dimensional family
since $\dim C\cup \cup_i \beta_i=1$ and
if one point $t$ of $C\cup \cup_i \beta_i$ is fixed,
then such pairs of conics such that
$t\in \overline{r}_1\cap \overline{r}_2$ 
form a one-dimensional family.
For a general pair of $r_1$ and $r_2$,
the number of the sets of $n$ mutually conjugate conics
including $r_1$ and $r_2$ is finite since
$D_{r_1}$ and $D_{r_2}$ has no common component.
Thus $\{q_i\}$ does not contain such a pair by generality
whence $\{q_i\}$ satisfies (3).   
\end{proof}

\section{Relation with Mukai's result}
\label{section:v22}
Here we sketch how the argument goes on if $d=5$ and
explain a relation of our result with
Theorem \ref{v22}.

Assume that $d=5$.
Associated to the birational morphism
$\sH_2\rightarrow\check{\mP}^{2}$,
there exists
a non finite birational morphism  
\[
\Phi\colon \widetilde{A}\rightarrow A_{22}:=\VSP(F_4,6)
\subset \Hilb^{6} \check{\mP}^{2},
\] 
which fits into the following diagram: 
\begin{equation*}
\xymatrix{ & & \widetilde{A} \ar[dl]_{\rho} \ar[dr]^{\rho'}
\ar @/^2pc/ [ddrr]^{\Phi} & &\\
& A¥ \ar[dl]_{f}  &\dashrightarrow & 
 A' \ar[dr]^{f'} & \\
 B &  &  & & A_{22}, }
\end{equation*}

\noindent 
where 
\begin{itemize}
\item
$A_{22}$ is a smooth prime Fano threefold of genus twelve,
\item
$\rho'$ is the blow-down of the three $\rho$-exceptional divisors 
$E_i$ ($i=1,2,3$)
over the strict transforms $\beta'_i$ in the other direction.
In other words, $A\dashrightarrow A'$ is the flops of $\beta'_1$, $\beta'_2$
and $\beta'_3$ (cf. Lemma \ref{sullebisecanti}), and
%\item
%$V_{22}$ is the smooth prime Fano $3$-fold of genus $12$, and
\item
the morphism $f'$ 
%is given
%by the strict transform of the linear system $|3f^*S-2E|$ 
contracts the strict transform of 
the unique hyperplane section $S$ containing $C$
(see Corollary \ref{cor:hull})
to a general line on $A_{22}$.
\end{itemize}
The rational map $A_{22}\dashrightarrow B$ is the famous double projection
of $A_{22}$ from a general line $m$ 
first discovered by Iskovskih (see \cite{I2}).

We explain how $f'$ and $\rho'$ are interpreted in our context. 
As we remarked after the proof of Theorem \ref{thm:H_2},
the morphism $\sH_2\to \check{\mP}^2$ defined by $|D_l|$ contracts
three curves $D_{e_i}$ which parameterize conics
intersecting $\beta'_i$.
By noting $S$ is covered by the images of such conics,
this corresponds to that 
the morphism $f'$ contracts the strict transform of $S$.

We can see that any conic on $A$ except one belonging to $D_{e_i}$
corresponds to that on $A_{22}$ in the usual sense, and
the component of Hilbert scheme of $A_{22}$ parameterizing
conics is naturally isomorphic to $\check{\mP}^2$.
The three conics on $A_{22}$ corresponding to the images of $D_{e_i}$
are $\beta''_i\cup m$, where $\beta''_i$ are the images of the flopped curve
corresponding to $\beta'_i$.

Let $a\in E_i$. Then six conics on $A$ attached to $a$
are $\xi_{ij1}$ ($j=1,2$), a conic $q_a$ from $D_{e_i}$ and
three conics from $e_i$ 
(see the remark at the end of \ref{subsection:finite}).
Moreover, if $a$ moves in a fiber $\gamma$ 
of the other projection $E_i\to \mP^1$,
then only the conic $q_a$ from $D_{e_i}$ varies.
By the contraction $\sH_2\to \check{\mP}^2$, 
there is no difference among points on $\gamma$.
This is the meaning of the contraction $\rho'$ of $E_i$
in the other direction.

Finally we remark that $\sH_1$ is also naturally isomorphic to    
the component of Hilbert scheme of $A_{22}$ parameterizing lines.

\section{Appendix}
\label{section:app}
We give a quick review of basic facts on the theory of polarity.
The main references are \cite[\S 1 and \S 2]{DK} and \cite[\S 2]{doldual}.
\begin{itemize}
\item
Denote by $\mathrm{Sym}_m V$
the image of the linear map
\begin{eqnarray*}
\check{V}^{\otimes m}&\to& \check{V}^{\otimes m}\\
t&\mapsto&  \sum_{\sigma\in S_m} \sigma(t).
\end{eqnarray*}
The map $\check{V}^{\otimes m}\to \mathrm{Sym}_m V$
is decomposed as
$\check{V}^{\otimes m}\stackrel{s_m}{\to} S^m \check{V}\stackrel{p_m}{\to} \mathrm{Sym}_m V$, where $s_m$ is the natural quotient map.
Denote by $\mathrm{pl}_m\colon S^m \check{V}{\to} \mathrm{Sym}_m V$
the map obtained from $p_m$ by dividing $m!$.
This is called {\em{the polarization map}}.
Let $r_m\colon \mathrm{Sym}_m V\hookrightarrow
\check{V}^{\otimes m}\stackrel{s_m}{\to} S^m \check{V}$
be the natural map.
Then it holds that
$\mathrm{pl}_m\circ r_m=r_m\circ \mathrm{pl}_m=\mathrm{id}$.
\item
For $F\in S^m \check{V}$, set $\widetilde{F}:=\mathrm{pl}_m(F)$.
Then $F(x)=\widetilde{F}(x,x,\dots,x)$ for $x\in V$.
\item   
For $F\in S^m \check{V}$ and $a\in V$,
set
$P_a(F)(x):=\widetilde{F}(a,x,\dots,x)$.
It is easy to varify
\[
P_a(F)=
\frac{1}{m}\sum_i a_i \frac{\partial F}{\partial {x}_i},
\]
where $a_i$ and $x_i$ are coordinates of $a$, and
on $V$ respectively.
Similarly, by setting
$P_{a,b,\dots,c}(F):=
\widetilde{F}(a,b,\dots,c,x,\dots,x)$
(the number of $a,b,\dots,c$ is $k$),
it holds
\[
P_{a,b,\dots,c,x,\dots,x}(F)=
\frac{(m-k)!}{m!}\sum_{i_1,\dots,i_k} a_{i_1}b_{i_2}\cdots c_{i_k} 
\frac{\partial^k F}{\partial {x}_{i_1}\cdots \partial {x}_{i_k}}.
\]
This is called {\em{the mixed polar of $F$ with respect to $a,b,\dots,c$.}}

It is possible to regard this as
the pairing between $F\in S^m \check{V}$ and $ab\cdots c\in S^k V$.
By extending this pairing,
we have
\begin{eqnarray*}
S^k {V}\times S^m \check{V}&\to& S^{m-k} \check{V}\\
(G,F)&\mapsto& P_G(F).
\end{eqnarray*}
Further, by fixing $F$, we can write  
\begin{eqnarray*}
\ap^k_{F}\colon S^k V& \to & S^{m-k} \check{V}\\
G & \mapsto & 
P_G(F).
\end{eqnarray*}
This is called the {\em{apolarity map}}.

When $m=k$, this pairing is sometimes denoted by
$\langle G,F\rangle$ and is called
the {\em{apolarity pairing}}.
\item
The following is a basic property of the apolarity pairing:
\[
\langle F, ab\cdots c\rangle=\widetilde{F} (a,b,\dots,c),
\]
where the number of $a,b,\dots,c$ is $m$.
In particular,
\[
\langle F, a^m\rangle=\widetilde{F} (a,a,\dots,a)=F(a).
\]
\item
If $m=2k$, then 
$F$ is said to be {\em{non-degenerate}} 
if
\[
\ap^k_{F}\colon
S^k V\to S^k \check{V}
\]
is an isomorphism. 
In this case, 
there is $\check{F}\in S^k V$
such that 
\[
{\ap^k_F}^{-1}=\ap^k_{\check{F}}.
\]
$\check{F}$ is called {\em{the form dual to ${F}$}}.
\item
Usually, we consider the apolarity maps in the projective setting.
Namely, we consider $a\in \mP_*V$ rather than $a\in V$, etc.
In this situation, we denote by $H_a\in V$ an element corresponding
to $a\in \mP_*V$, which is unique up to scalar.
By abuse of notation, we sometimes continue to write
$P_a(F)$ rather than $P_{H_a}(F)$.
\end{itemize}

\end{document}